\documentclass[12pt, reqno]{amsart}

\usepackage{graphics}

\newtheorem{theorem}{Theorem}
\newtheorem{proposition}[theorem]{Proposition}
\newtheorem{lemma}[theorem]{Lemma}
\newtheorem{corollary}[theorem]{Corollary}

\theoremstyle{definition}
\newtheorem{definition}[theorem]{Definition}
\newtheorem{example}[theorem]{Example}

\theoremstyle{remark}

\numberwithin{equation}{section}

\newcommand{\Z}{{\mathbf Z}}
\newcommand{\Q}{{\mathbf Q}}
\newcommand{\R}{{\mathbf R}}

\newcommand{\kk}{{\mathbf k}}
\newcommand{\cat}{{\rm {cat }}}

\newcommand{\ind}{{\mbox{\rm ind}}}
\newcommand{\RP}{\mathbf {RP}}
\newcommand{\const}{\rm {const}}
\newcommand{\cl}{\rm {cl}}
\begin{document}

\title{Topology of billiard problems, II}         
\author{Michael Farber}        
\address{Department of Mathematics, Tel Aviv University, Tel Aviv, 69978, Israel}
\email{farber@math.tau.ac.il}

\date{\today} 
         
\begin{abstract} In this paper we 
give topological lower bounds on the numbers of periodic and closed trajectories in strictly convex smooth billiards
$T\subset \R^{m+1}$. Namely, for given $n$ we estimate the number of $n$-periodic billiard trajectories in $T$, and also the
number of billiard trajectories,
which start and end at a given point $A\in \partial T$ 
and make a prescribed number $n$ of reflections at the boundary $\partial T$ of the billiard domain.
We use variational reduction, admitting a finite group symmetries, and apply topological approach based on 
equivariant Morse and Lusternik - Schnirelman theories.
\end{abstract}
\thanks{Partially supported by the US - Israel Binational Science Foundation and by the Minkowski Center for Geometry}

\maketitle

\section{\bf Introduction}     

Let $X\subset \R^{m+1}$ be a closed smooth strictly convex hypersurface. 
We will consider the billiard system in the $(m+1)$-dimensional convex body $T$, bounded by $X$. 
Recall, that we view the billiard ball as a point, which moves in $T$ in a straight line, except when it hits $X= \partial T$,
where it rebounds, making the angle of incidence equal the angle of reflection.

G.D. Birkhoff \cite{Bi} studied periodic billiard trajectories in plane convex billiards.
Papers \cite{Ba}, \cite{FT} deal with the problem of estimating the number of periodic 
trajectories in convex billiards in $\R^{m+1}$, where $m>1$. 
In \cite{F} we studied the number of billiard trajectories
having fixed distinct end points and making a prescribed number of reflections.

The purpose of this paper, which continues \cite{F}, is twofold.
Firstly, we obtain estimates on the number of closed billiard trajectories,
which start and end at a given point $A\in X$ and
make a prescribed number $n$ of reflections at the hypersurface $X$. This problem may look as a special case
of the fixed end billiard problem \cite{F}, but as we show, presence of symmetry allows to get much stronger estimates
than \cite{F}.
Secondly, we give a linear in $n$
estimate on the number of $n$-periodic trajectories. 

The following Theorem \ref{thm2} gives an estimate on the number of closed billiard trajectories. 
It deals with $\Z_2$-orbits
of billiard trajectories. Any such $\Z_2$-orbit is determined by a sequence of 
reflection points $x_1, x_2, \dots, x_n\in X$, such that $x_i\ne x_{i+1}$ for 
$i=1, \dots, n-1$ and $x_1\ne A$, $x_n\ne A$. The reverse sequence $x_n, x_{n-1}, \dots, x_1$ determines the same $\Z_2$-orbit. 

\begin{theorem}\label{thm2}
Let $X\subset \R^{m+1}$ be a closed smooth strictly convex hypersurface, $A\in X$. 

(I) For any $n\ge 1$ the number of distinct $\Z_2$-orbits of closed billiard trajectories inside $X$, 
which start and end at $A$, and make $n$ reflections is at least 
\begin{eqnarray}
\begin{array}{ll}
n,            &\mbox{if $m\ge 3$ is odd}, \\ \\
{}[n/2] +1, &\mbox{if $m\ge 2$ is even}. 
\end{array}
\end{eqnarray}

(II) For any even $n\ge 2$ the number of distinct $\Z_2$-orbits of closed billiard trajectories inside $X$, 
which start and end at $A$, and make $n$ reflections is at least 
\begin{eqnarray}
\begin{array}{lll}
 [\log_2 n] +m-1,&\mbox{if $m\ge 3$ is odd},\\ \\
 {}[\log_2n]+m-2, &\mbox{if $m\ge 2$ is even and $n\ge 4$},\\ \\
m, & \mbox{if $m\ge 2$ is even and $n=2$.}
\end{array}
\end{eqnarray}

(III) If  $n\ge 2$ is even, and the billiard data $(X, A, n)$ is generic (cf. below), then 
the number of distinct $\Z_2$-orbits of closed billiard trajectories inside $X$, which start and end at $A$, 
and make
$n$ reflections is at least 
\begin{eqnarray}
 mn/2.
\end{eqnarray}
\end{theorem}

First we explain the {\it genericity} assumption in statement {\it (III)}. The billiard data $(X,A,n)$ determines a continuous 
function
\begin{eqnarray}\label{func}
X^{\times n}\to \R, \quad (x_1, \dots, x_n)\mapsto \sum\limits_{i=0}^n |x_i-x_{i+1}|,
\end{eqnarray}
(the total length),
where we understand $x_0=A=x_{n+1}$. This function is smooth at all configurations $(x_1, \dots, x_n)\in X^{\times n}$ with 
$x_i\ne x_{i+1}$ for $i=0, \dots, n$. The data $(X,A,n)$ is {\it generic} if any critical configuration 
$(x_1, \dots, x_n)\in X^{\times n}$ of the total length function (\ref{func}), satisfying the above condition $x_i\ne x_{i+1}$, 
is Morse. 
Compare \cite{Ba}, \cite{FT}. 

Statements (I) and (II) give different lower bounds on the number of closed billiard trajectories. (I) is linear in $n$;
it is better than (II) for large $n$. On the other hand, (II) may be better than (I)
if the dimension $m=\dim X$ of the boundary of billiard domain is large. 

Let us compare Theorem \ref{thm2} with the lower bound on the number of billiard trajectories with fixed distinct end points, 
obtained in \cite{F}. In Theorem \ref{thm2} we speak about $\Z_2$-orbits of billiard trajectories. 
Each $\Z_2$-orbit contains one or two billiard trajectories. For $n$ even, each $\Z_2$-orbit contains precisely
two distinct billiard trajectories. Hence, we see that for $n$ even statement (I) of Theorem \ref{thm2}
predicts twice number of closed billiard trajectories, compared to the
estimate of \cite{F} for the billiard trajectories with fixed ends. 
Also, for large $m$ statements (II) and (III) give much larger lower bounds than corresponding estimates of Theorem 1 of 
\cite{F}.

Statement (III) includes the case $m=1$ (the plane billiards) and gives the estimate $n/2$.
The billiard in the unit circle has precisely $n/2$ orbits of closed billiard trajectories with a given initial point.

It is reasonable to expect that for any even $n\ge 2$ 
the number of distinct $\Z_2$-orbits of closed billiard trajectories inside $X$, 
which start and end at $A$, and make $n$ reflections is at least 
\begin{eqnarray}
\begin{array}{lll}
n +m-1,&\mbox{if $m\ge 3$ is odd},\\ \\
 {}n/2+m-1, &\mbox{if $m\ge 2$ is even},\\
\end{array}
\end{eqnarray}
Such estimate would imply both statements (I) and (II) of Theorem \ref{thm2}. 
The methods of this paper do not to prove this assertion, although the gap looks very small. 

The proof of Theorem \ref{thm2} is based on a computation of the cohomology ring of a relevant configuration space
of points on the sphere $S^m$.
We apply the technique of the critical point theory, based on the cup-length estimates together with a refinement, 
suggested by E. Fadell and S. Huseini \cite{FH},
related to the notion of category weight of cohomology classes.   

Next we state the main result concerning $n$-periodic trajectories.

\begin{theorem}\label{thm5}
Let $X\subset \R^{m+1}$ be a smooth strictly convex hypersurface.
For any odd prime $n$ there exist at least 
\begin{eqnarray}
\begin{array}{ll}
n, &\mbox{if $m$ is odd},\\
(n+1)/2,&\mbox{if $m$ is even}
\end{array}
\end{eqnarray}
distinct $D_n$-orbits of $n$-periodic billiard trajectories inside $X$.
\end{theorem}

This theorem complements the results of \cite{FT}.
In \cite{FT} it is shown that for $m\ge 3$ and odd $n$ the number of 
distinct $D_n$-orbits of $n$-periodic 
billiard trajectories inside $X\subset \R^{m+1}$ is not less than $[\log_2(n-1)]+m$ and it at least 
$(n-1)m$
for generic billiards $X\subset \R^{m+1}$. These results from \cite{FT}
are similar to statements (II) and (III) of Theorem \ref{thm2}.
Theorem \ref{thm5} above has several advantages compared to \cite{FT}. 
It gives a linear in $n$ estimate, which is better for large $n$
than the logarithmic estimate of \cite{FT}. Also, it allows the case $m=2$ which corresponds to convex 
billiards in 3-dimensional Euclidean space. On the other hand, the result of \cite{FT} is better for large $m$. 

The proof of Theorem \ref{thm5} is based on a computation of the cohomology rings of cyclic configuration spaces
of spheres with rational coefficients. The case of $\Z_2$-coefficients was computed in \cite{FT}.

I would like to thank S. Tabachnikov for useful discussions.

\section{\bf Cohomology of the closed string configuration spaces of spheres}

Let 
\begin{eqnarray}\label{closedstring}
G_n =G(S^m;A,A,n)
\end{eqnarray}
denote the closed string configuration space of $S^m$, i.e. the space of all 
configurations $(x_1, \dots, x_n)$, where $x_i\in S^m$, such that $x_1\ne A$, $x_n\ne A$ and 
$x_i \ne x_{i+1}$ for all $i=1, \dots, n-1$. There is a natural involution
\begin{eqnarray}\label{reflection}
T: G_n \to G_n, \quad T(x_1, \dots, x_n) = (x_n. x_{n-1}, \dots, x_1),
\end{eqnarray}
which will be important for the sequel.

\begin{theorem}\label{openstring5}
The cohomology group $H^i(G_n;\Z)$ is nonzero only in dimensions 
$$i=0,\,  (m-1),\,  2(m-1), \dots, (n-1)(m-1)$$ and for these
values $i$ the group $H^i(G_n;\Z)$ is free abelian of rank 1. 
One may choose additive generators 
\[\sigma_i \in \, H^{i(m-1)}(G_n;\Z), \quad i=0, 1, \dots , n-1,\]
such that for $m \ge 3$ odd, the multiplication is given by 
\begin{eqnarray}\label{prod1}
\sigma_i\sigma_j \, =\, \left \{
\begin{array}{l}
\displaystyle{\frac{(i+j)!}{i!\cdot j!}}\cdot \sigma_{i+j},\quad\mbox{if}\quad i+j\le n-1,\\ \\
0, \quad \mbox{if}\quad  i+j>n-1 
\end{array}
\right .
\end{eqnarray}
and for $m\ge 2$ even, it is given by
\begin{eqnarray}\label{prod2}
 \sigma_i\sigma_j \, =\, \left \{
\begin{array}{l}
\displaystyle{\frac{[(i+j)/2]!}{[i/2]!\cdot [j/2]!}}\cdot \sigma_{i+j},\quad\mbox{if $i+j\le n-1$ and $i$ or $j$ is even,}\\ \\
0, \quad \mbox{if either $i+j>n-1$, or both $i$ and $j$ are odd.}
\end{array}
\right.
\end{eqnarray}
Reflection (\ref{reflection}) acts for $m> 1$ odd by
\begin{eqnarray}\label{reflection1}
T^\ast(\sigma_i) = (-1)^i \sigma_i, 
\end{eqnarray}
and for $m>1$ even by
\begin{eqnarray}\label{reflection2}
 T^\ast(\sigma_{i}) = (-1)^{[i/2]+ni}\sigma_i,\quad i=0, 1, \dots, n-1.
\end{eqnarray}
\end{theorem}

\begin{proof} Consider the map
$$G_n=G(S^m;A,A,n)\,  \to\,  S^m - A, \quad (x_1, \dots, x_n) \mapsto x_n.$$
It is a smooth fibration with fiber $G(S^m;A,B,n-1)$, where $A\ne B$. 
Since the base $S^m-C$ is contractible,
we obtain that the inclusion
\begin{eqnarray}\label{equivalence}
G(S^m;A,B, n-1) \subset G_n 
\end{eqnarray}
is a homotopy equivalence. Hence, the integral cohomology ring of $G_n$ coincides with
$H^\ast(G(S^m;A,B,n-1);\Z)$, which we calculate below. 

Theorem 8 from \cite{F} describes algebra $H^\ast(G(S^m;A,B,n-1);\kk)$, where
$\kk$ is an arbitrary field. From this description it is clear that the dimension of the cohomology does not depend on field $\kk$.
Therefore, we conclude that the integral cohomology $H^i(G(S^m;A,B,n-1);\Z)$ has no torsion;
it is a free abelian group of rank one for $i=r(m-1)$, where $r=0, 1,\dots, n-1$, and it vanishes for all other values 
of $i$.

Let $C\in S^m$ be a point distinct from $A$ and $B$. We obtain an inclusion of configuration spaces
$\phi^\ast : G(S^m-C; A, B, n-1) \to G(S^m; A, B,n-1)$, where we identify $S^m-C$ with $\R^m$.
The cohomology algebra $H^\ast(G(\R^m;A,B,n-1);\Z)$ has generators $s_0, \dots, s_{n-1}$ and the full list of relations
was described in Proposition 7 of \cite{F}.
From remark 9 in \cite{F} we know that the induced map $\phi^\ast$ on cohomology with arbitrary field of coefficients
$\kk$ is injective. This implies that the induced map on integral cohomology
\begin{eqnarray}
\phi^\ast: H^\ast(G(S^m;A,B,n-1);\Z) \to H^\ast(G(\R^m;A,B,n-1);\Z)
\end{eqnarray}
is injective and $\phi^\ast$ maps indivisible classes from $H^\ast(G(S^m;A,B,n-1);\Z)$ into indivisible classes in 
$H^\ast(G(\R^m;A,B,n-1);\Z)$. 

We claim that for any $r=0, 1, \dots, n-1$ there exists an indivisible class
$$\sigma_r\in H^{r(m-1)}(G(S^m;A,B,n-1);\Z),$$
such that 
\begin{eqnarray}\label{multiple}
 \phi^\ast(\sigma_r) = \left\{
\begin{array}{ll}
\sum\limits_{0\le i_1<\dots < i_r < n}s_{i_1}\dots s_{i_r} & \mbox{for $m$ odd},\\ \\
(-1)^{[r/2]+nr}\cdot\sum\limits_{0\le i_1<\dots < i_r < n}(-1)^{i_1+\dots +i_r} s_{i_1}\dots s_{i_r} & \mbox{for $m$ even}
\end{array}
\right.
\end{eqnarray}
(compare with formulae (4.3) and (4.4) from \cite{F}).
 Indeed, applying Remark 9 from \cite{F} with $\kk=\Q$, we see that the image of the generator of the group 
$H^{r(m-1)}(G(S^m;A,B,n-1);\Z)\simeq \Z$ under homomorphism $\phi^\ast$ equals an integral 
multiple of the expression in the RHS of (\ref{multiple}). Since the classes in RHS of (\ref{multiple}) are indivisible,
and since we know that $\phi^\ast$ maps indivisible classes to indivisible classes,
we conclude that a there exists generator $\sigma_r$ with the required property.

The product formulae (\ref{prod1}) and (\ref{prod2}) for classes $\sigma_r$ follow since they hold for the 
products $\phi^\ast(\sigma_i)\phi(\sigma_j) \in H^\ast(G(\R^m;A,B,n-1);\Z)$ as can be easily checked using the arguments of the
proof of Theorem 8 from \cite{F}.

Now we want to find the action of the reflection $T: G_n \to G_n$ on classes $\sigma_i$. 
It is clear that $T^\ast(\sigma_i) = \pm \, \sigma_i$, and we need to calculate the sign. 
Consider the following diagram of natural inclusions
$$
\begin{array}{ccc}
G(\R^m;A,B,n-1) & \to & G(S^m;A,B,n-1)\\ 
\downarrow & & \downarrow \\ 
G(\R^m;A,A,n) & \to & G(S^m;A,A,n)
\end{array}
$$
(where $\R^m=S^m-C$ as above) and the induced diagram of cohomology groups
$$
\begin{array}{ccc}
 H^\ast(G(\R^m;A,A,n);\Z)&\stackrel{\gamma} \leftarrow &H^\ast(G(S^m;A,A,n);\Z) \\ 
\beta \downarrow & &\alpha \downarrow \simeq \\ 
 H^\ast(G(\R^m;A,B,n-1);\Z)& \stackrel {\phi^\ast}\leftarrow & H^\ast(G(S^m;A,B,n-1);\Z)
\end{array}
$$
$\alpha$ is an isomorphism and $\phi^\ast$ is injective. To understand $\beta$, note that 
$G(\R^m;A,A,n)$ is homotopy equivalent to the cyclic configuration space $G(\R^m, n+1)$ (cf. \cite{FT})
and so the cohomology  $H^\ast(G(\R^m;A,A,n);\Z)$ has $(m-1)$-dimensional generators $s_0, s_1, \dots, s_n$
which satisfy relations of Proposition 2.2 from \cite{FT} (we shift indices for convenience). Proof of Proposition 7 
in \cite{F} shows that $\beta(s_i)=s_i$ for $i=0, 1, \dots, n-1$ and $\beta(s_n)=0$. Hence, $\beta$ is an epimorphism
with kernel equal the ideal generated by $s_n$. 

The reflection $T$ acts also on $G(\R^m;A,A,n)$ (by formula (\ref{reflection}). 
It is clear that the induced map $T^\ast: H^\ast(G(\R^m;A,A,n);\Z)\to  H^\ast(G(\R^m;A,A,n);\Z)$ acts on the generators
$s_i$ as follows 
\begin{eqnarray}\label{reflection3}
T^\ast(s_i) = (-1)^m s_{n-i}, \quad\mbox{where}\quad  i=0, 1, \dots, n.
\end{eqnarray}

Now we may calculate $T^\ast(\sigma_r)$, where $r=1, \dots, n-1$. 
Fix a subsequence $0<i_1<\dots <i_r<n$ (we avoid indices $0$ and $n$).
Suppose first that $m$ is odd. Then $\phi^\ast(\alpha(\sigma_r))$ contains monomial $s_{i_1}s_{i_2}\dots s_{i_r}$; 
therefore $\gamma(\sigma_r)$ contains the same monomial with coefficient 1. Then $T^\ast(\gamma(\sigma_r))$ contains
monomial $s_{n-i_r}s_{n-i_{r-1}}\dots s_{n-i_1}$ with coefficient $(-1)^{mr}=(-1)^r$. The last monomial appear in 
$\gamma(\sigma_r)$ with coefficient 1. Since we know that $T^\ast(\sigma_r) = \pm \, \sigma_r$ we conclude that 
$T^\ast(\sigma_r) = (-1)^r  \sigma_r$.

Assume now that $m$ is even. 
Then 
$\phi^\ast(\alpha(\sigma_r))$ contains monomial $s_{i_1}s_{i_2}\dots s_{i_r}$ with coefficient 
$$(-1)^{[r/2]+nr+i_1+\dots+ i_r}.$$ 
Applying $T^\ast$ and using (\ref{reflection3}) we see that the monomial $s_{n-i_r}s_{n-i_{r-1}}\dots s_{n-i_1}$ 
appears in $T^\ast(\gamma(\sigma_r))$ with coefficient 
$$(-1)^{nr +i_1+\dots+ i_r}$$ 
and in 
$\gamma(\sigma_r)$ with coefficient 
$$(-1)^{[r/2]+i_1+\dots+ i_r}.$$
This shows that 
$T^\ast(\sigma_r) = (-1)^{[r/2]+nr}  \sigma_r$.
\end{proof}

\section{\bf Calculation of equivariant cohomology}

Our purpose in this section is to compute the cohomology of $G_n/\Z_2$,
the factor-space 
of the space of closed string configurations $G_n=G(S^m;A,A,n)$
with respect to $\Z_2$-action given by the reflection $T: G_n\to G_n$.
For $n$ even $T$ acts freely, and $H^\ast(G_n/\Z_2;\Z)$ coincides with the
equivariant cohomology of $G_n$.

 The problem is trivial for $m=1$; 
therefore everywhere in this section we will assume that $m>1$.

To compute the equivariant cohomology we will apply the method of the Morse theory. 
Namely, we will consider the simplest billiard
in the standard unit sphere $S^m\subset \R^{m+1}$ and the function of negative total length
\begin{eqnarray}\label{function}
L: G_n=G(S^m;A,A,n) \to \R, \quad (x_1, \dots, x_n)\mapsto -\sum\limits_{i=0}^n |x_i-x_{i+1}|.
\end{eqnarray}
Here we understand $x_0=x_{n+1}=A$.
The critical points of $L$ are the billiard trajectories in $S^m$, which start and end at $A$ and make $n$ reflections.
All such trajectories can easily be described. 

Namely, fix a vector $a\in S^m$, $a\perp A$, orthogonal to $A$ and an angle
\begin{eqnarray}\label{angle}
\psi_k = \frac{2\pi k}{n+1}, \quad k=1, 2, \dots, [(n+1)/2].
\end{eqnarray}
This choice $(a, \psi_k)$ determines the following billiard trajectory $(x_1, \dots, x_n)$, where
\[x_j = A\cos(j\psi_k) +a\sin(j\psi_k), \quad j=1, \dots, n.\]
Note that for $n$ odd the trajectory determined by the pair $(a, \psi_{(n+1)/2})$ does not depend on $a$;
it has the form $(x_1, \dots, x_n)$,
where $x_j =A$ for $j$ even and $x_j=-A$ for $j$ odd. 

We will denote by 
$$V_p\subset G_n, \quad p=0, 1, \dots, [(n-1)/2]$$
 the variety of trajectories determined by all pairs $(a, \psi_{k})$, where
\begin{eqnarray}\label{conven}
k=[(n+1)/2]-p
\end{eqnarray}
 and $a\perp A$ is an arbitrary point of the sphere $S^{m-1}\subset S^m$ orthogonal to $A$. 

If $n$ is even then every submanifold $V_p$ is diffeomorphic to sphere $S^{m-1}$. 

If $n$ is odd, then $V_0$ is a single point and 
$V_1, \dots, V_{[(n-1)/2]}$ are diffeomorphic to the sphere $S^{m-1}$.

The following statement is similar to Proposition 3.1 of I.K. Babenko \cite{Ba}.

\begin{proposition}\label{prop1}
Each $V_p\subset G_n$ is a nondegenerate critical submanifold of function $L$ in the sense of Bott. 

If $n$ is even then index of each $V_p$ equals $2p(m-1)$ for $p=0, 1, \dots, (n-2)/2$. 

If $n$ is odd, then index of $V_0$ equals $0$ and for $p=1, \dots, (n-1)/2$
the index
of $V_p$ equals $(2p-1)(m-1)$.
\end{proposition}
\begin{proof} Let $e_1, \dots, e_{m+1}\in \R^{m+1}$ be an orthonormal base. We may assume that $A=e_1$.
We want to calculate the Hessian  of function $L$ at a billiard trajectory $c_k= (x_1, \dots, x_n)\in G_n$, where
\[x_j =\cos(\psi_k) e_1 +\sin (\psi_k)e_2, \quad j=1, \dots, n, \]
and 
\begin{eqnarray*}
\psi_k = \frac{2\pi k}{n+1}, \quad k=1, \dots, [(n+1)/2].
\end{eqnarray*}
Let $x_j^\perp$ denote the orthogonal vector to $x_j$ lying in the $(e_1, e_2)$-plane, i.e. 
\[x_j^\perp = \cos(\psi_k+\pi/2) e_1 +\sin (\psi_k+\pi/2)e_2.\]
Any tangent vector $Y\in T_{c_k}G_n = \oplus_j T_{x_j}S^m$ is determined by numbers $\mu_{r, j}\in \R$,
where $r=0, 1, \dots, m-1$ and $j=1, \dots, n$, 
such that the component of $Y$ in $T_{x_j}S^m$ equals 
\[\mu_{0,j}x_j^\perp + \sum\limits_{r=1}^{m-1} \mu_{r,j}e_{r+2}.\]
A direct calculation of Hessian $H(L)_{c_k}(Y, Y)$ of $L$ gives the following quadratic form in variables $\mu_{r,j}$:
\begin{eqnarray}\label{hessian}
\begin{array}{lll}
H(L)_{c_k}(Y,Y) &=& \displaystyle{\frac{1}{2} \sin(\psi_k/2)} \cdot \sum\limits_{j=0}^n (\mu_{0,j }-\mu_{0, j+1})^2 + \\ \\
&+& \displaystyle{\left({2\sin(\psi_k/2)}\right)^{-1}} \cdot\left[\sum\limits_{r=1}^{m-1} Q_{\psi_k}(\mu_{r,1}, \dots, \mu_{r, n})\right],
\end{array}
\end{eqnarray}
where in the first sum we understand $\mu_{0,0}=0=\mu_{0,n+1}$ and in the second sum the symbol
$Q_\psi(y_1, \dots, y_n)$ denotes the following quadratic form
\[Q_\psi(y_1, \dots, y_n) = -2\cos(\psi)\cdot\sum\limits_{i=1}^n y_i^2 + 2\sum\limits_{i=1}^{n-1}y_iy_{i+1}.\]

We see that the Hessian splits as a direct sum of $m$ quadratic forms corresponding to different values $r=0, 1, \dots, m-1$.
The terms involving $\mu_{0,j}$ (the first sum) give a positive definite quadratic form. 
The rest $m-1$ form are identical and their index and nullity equal the index and nullity of $Q_{\psi_k}$.
Hence we conclude that 
the index and nullity of the Hessian
equals $m-1$ times the index and nullity of the form $Q_{\psi_k}$. 

In order to calculate the index of $Q_{\psi_k}$ we observe that 
the eigenvalues of the following symmetric $n\times n$-matrix 
$$
 \left[
\begin{array}{cccccc}
0 & 1& 0 & \dots & 0 & 0\\
1 & 0& 1 & \dots & 0 & 0\\
0 & 1& 0 & \dots & 0 & 0\\
\dots & \dots & \dots & \dots & \dots &\dots \\
0 & 0 & 0 & \dots & 0& 1\\
0 & 0& 0 & \dots & 1 & 0
\end{array}
\right]
$$
are given by
\begin{eqnarray}\label{values}
\lambda_s = 2\cos\left(\frac{\pi s}{n+1}\right), \quad s =1, 2, \dots, n
\end{eqnarray}
and the eigenvector $(v_{1,s}, \dots, v_{n,s})$ corresponding to $\lambda_s$ is given by 
\begin{eqnarray}\label{vector}
v_{j,s} = \sin\left(\frac{\pi j s}{n+1}\right), \quad j=1, \dots, n.
\end{eqnarray}
This claim can be checked directly. 

Therefore the eigenvalues of $Q_{\psi_k}$ are 
\begin{eqnarray}\label{values1}
2\left[\cos\left(\frac{\pi s}{n+1}\right)-\cos\left(\frac{2\pi k}{n+1}\right)\right], \quad s =1, 2, \dots, n
\end{eqnarray}
and the eigenvectors of $Q_{\psi_k}$ are given by (\ref{vector}).

Hence the index of $Q_{\psi_k}$ equals the number of integers $s$, such that $2k<s\le n$, which is $n-2k$ for
$2k\le n$ and $0$ if $k=(n+1)/2$ and $n$ is odd. 
Since (according to (\ref{conven})), $k=[(n+1)/2]-p$, we conclude that index of $Q_{\psi_k}$ equals 
\[n-2k =n-2([(n+1)/2]-p) = \left\{
\begin{array}{ll}
2p, &\mbox{if $n$ is even}\\
2p-1, &\mbox{if $n$ is odd}.
\end{array}
\right.\]
The special case $k=(n+1)/2$ for $n$ odd corresponds to $p=0$; in this case the index and nullity 
of $Q_{\psi_k}$ equal $0$.

From (\ref{values1}) we see that the nullity of $Q_{\psi_k}$ equals 1 for any $k$ unless $n$ is odd and $k=(n+1)/2$.

The discussion above proves that on any critical submanifold $V_p$ 
the dimension of the kernel of Hessian of $L$ equals the dimension of $V_p$; hence all submanifolds
$V_p$ are nondegenerate in the sense of Bott and their indices are as stated.
\end{proof}

The normal bundle $\nu(V_p)$ splits as a direct sum $\nu_+(V_p)\oplus \nu_-(V_p)$ of the positive and negative normal bundles with respect to the Hessian of $L$. 
One may describe the negative normal bundle $\nu_-(V_p)$ as follows.

\begin{lemma} \label{lemma1}The negative normal bundle $\nu_-(V_p)$ to $V_p$ is
\begin{eqnarray}
\nu_-(V_p) \, =\, \left\{ 
\begin{array}{ll}
\underbrace{\xi\oplus \xi\oplus \dots \oplus\xi}_{2p\quad\mbox{times}},& \mbox{if $n$ is even},\\ \\
\underbrace{\xi\oplus \xi\oplus \dots \oplus\xi}_{2p-1\quad\mbox{times}},& \mbox{if $n$ is odd and $p>0$},
\end{array}
\right.
\end{eqnarray}
where $\xi$ denotes the tangent bundle of sphere $S^{m-1}$.
\end{lemma}
\begin{proof} Let $S^{m-1}\subset S^m$ be the equatorial sphere consisting of unit vectors orthogonal to $A$.
Any point $a\in S^{m-1}$ and an angle (\ref{angle}) determine a critical submanifold $V_p$. Fix an eigenvalue
$\lambda_s$ (given by (\ref{values})), such that expression (\ref{values1}) is negative. Consider the subbundle
$\nu_s(V_p)$ of the normal bundle $\nu(V_p)$ consisting of eigen vectors of the Hessian with eigenvalue $\lambda_s$. 
We want to show that $\nu_s(V_p)$ is isomorphic to $\xi$. This would clearly imply the Lemma.

Consider a billiard trajectory $c_k= (x_1, \dots, x_n)\in G_n$ in the plane of vectors $a$ and $A$, where
\[x_j =\cos(\psi_k) A +\sin (\psi_k)a, \quad j=1, \dots, n, \]
and 
\begin{eqnarray*}
\psi_k = \frac{2\pi k}{n+1}, \quad k= [(n+1)/2] -p.
\end{eqnarray*}
Denote by $\xi_a$ the $(m-1)$-dimensional subspace orthogonal to $a$ and $A$. 
We will show that there is an isomorphism between the fiber of $\nu_s(V_p)$ over $c_k$ and $\xi_a$,
which depends continuously on $a$. 

Let $v_j\in T_{x_j}S^m$, where $j=1, 2, \dots, n$,
be a sequence of tangent vectors. Using (\ref{hessian}) and (\ref{vector}) we find that a 
sequence of vectors $(v_1, \dots, v_n)$ belongs to the fiber of 
$\nu_s(V_p)$ over the configuration $c_k=(x_1, \dots, x_n)$ if and only if 
\begin{eqnarray}\label{first}
v_j\in \xi_a, \quad \mbox{and}\quad 
v_j = \frac{\sin\left(\displaystyle{\frac{\pi j s}{n+1}}\right)}{\sin\left(\displaystyle{\frac{\pi  s}{n+1}}\right)}\cdot v_1,
\quad j=1, \dots, n.\end{eqnarray}
We see that the first vector $v_1$ uniquely determines a tangent vector $(v_1, \dots, v_n)$ to a configuration $c_k$
in the eigen-deriction $\lambda_s$. Moreover, $v_1$ can be an arbitrary vector in $\xi_a$. 
\end{proof}

Since $\xi$ is orientable we obtain:

\begin{corollary}\label{orient}
 The negative normal bundle $\nu_-(V_p)$ is orientable.
\end{corollary}

Note that this Corollary is trivial for $m>2$ since then the sphere $S^{m-1}$ is simply connected.

\begin{corollary}\label{perfect1} The function $L: G_n\to \R$ (cf. (\ref{function})) is a perfect Bott function.
 \end{corollary}
\begin{proof} Note that the critical value $L(V_p)$ equals
$$L(V_p)\, =\, -2(n+1)\sin\left(\frac{2\pi k}{n+1}\right), \quad\mbox{where}\quad k=[(n+1)/2]-p.$$
Hence for $p<p'$ we have $L(V_p)<L(V_{p'})$.

Choose constants $c_0, c_1, \dots, c_{[(n-1)/2]}\in \R$ such that 
\[L(V_p)<c_p<L(V_{p+1}), \quad p=0, 1, \dots, [(n-3)/2],\quad\mbox{and}\quad L(V_{[(n-1)/2]})<c_{[(n-1)/2]}.\]
Each $F_p=L^{-1}(-\infty, c_p])\subset G_n$ is a compact manifold with boundary and we obtain a filtration
\[F_0\subset F_1\subset \dots \subset F_{[(n-1)/2]},\]
and the inclusion $F_{[(n-1)/2]}\to G_n$ is a homotopy equivalence (as follows easily from Proposition 4 in \cite{F}). 
Using Corollary \ref{orient} and the Thom isomorphism we obtain
\begin{eqnarray}
\begin{array}{lll}
H^j(F_p,F_{p-1};\Z) &\simeq & H^{j-\ind(V_p)}(V_p;\Z) \, =\\ \\
&=& \, \left\{
\begin{array}{ll}
\Z, & \mbox{if $j= \ind(V_p)$ or $j= \ind(V_p)+m-1$}\\ \\
0, & \mbox{otherwise}
\end{array}
\right.
\end{array}
\end{eqnarray}
This holds true also for $p=0$ if we understand $F_{-1}=\emptyset$.

Suppose that $n$ is even. Then cohomology group $H^j(F_p,F_{p-1};\Z)$ is isomorphic to $\Z$ for $j=2p(m-1)$ and for $j=(2p+1)(m-1)$
and vanishes for all other $j$. Comparing with additive structure of $H^\ast(G_n;\Z)$ given by Theorem \ref{openstring5}
we find
\begin{eqnarray}\label{perfect}
H^\ast(G_n;\Z) \simeq \bigoplus_{p=0}^{[(n-1)/2]} H^\ast(F_p,F_{p-1};\Z),
\end{eqnarray}
which means perfectness of $L$. 

Suppose now that $n$ is odd. Then $H^j(F_0,F_{-1};\Z)$ is $\Z$ for $j=0$ and vanishes for all other values of $j$.
If $p>0$ then 
$$
H^j(F_p,F_{p-1};\Z)\, \simeq \, \left\{
\begin{array}{l}
\Z,\quad\mbox{for $j=(2p-1)(m-1)$ or $j=2p(m-1)$}\\
0, \quad\mbox{otherwise}
\end{array}
\right.
$$
and thus perfectness (\ref{perfect}) also holds.

Alternatively, for $m>2$ perfectness (\ref{perfect}) follows without using Theorem \ref{openstring5}
by considering the spectral sequence of filtration $F_p$ 
$$E_1^{p,q} = H^{p+q}(F_p, F_{p-1};\Z) \Rightarrow H^{p+q}(G_n;\Z)$$
and observing that for any of its differential $d_r$, with $r\ge 1$, either the source or the target vanish.
Therefore, $E_1=E_\infty$. Moreover, every diagonal $p+q=c$ of $E_\infty$ contains at most one nonzero group.
If $m=2$ the differential $d_1$ has nonzero source and target, and so the above argument does not work. 
\end{proof}

From this point we will assume that $n$ is even. 

Then the reflection $T: G_n\to G_n$ acts freely and our purpose 
will be to calculate the cohomology of the factor-space $G_n'= G_n/\Z_2.$
Function (\ref{function}) is reflection invariant and so it determines a smooth function 
$$ L': G_n'\to \R.$$
The critical points of $L'$ form nondegenerate (in the sense of Bott) critical submanifolds
\[V'_0, V'_1, \dots, V'_{n/2-1},\]
where $V'_p = V_p/\Z_2$. Index of $V_p'$ equals $2p(m-1)$ (as follows from Proposition \ref{prop1}).
Since each $V_p$ can be identified with $S^{m-1}$ and under this identification
the reflection $T$ acts as the usual antipodal map, we see that each $V_p'$ is diffeomorphic to the projective space
$\RP^{m-1}$.  
\begin{corollary}\label{poincare}
The Poincar\'e polynomial of $G_n'=G(S^m;A,A,n)/\Z_2$ with coefficients in field $\Z_2$
is 
\[\frac{t^m-1}{t-1}\cdot \frac{t^{n(m-1)}-1}{t^{2(m-1)}-1},\]
and the sum of Betti numbers with coefficients in $\Z_2$ is $mn/2$.
\end{corollary}
\begin{proof} We will give here a simple proof working for $m>2$. The case $m=2$ will follow from Theorem \ref{thm3}
below.

Consider filtration $F_0\subset F_1\subset \dots \subset F_{n/2-1}\subset G_n$ as in the proof of Corollary \ref{perfect1}. Let $F_p'$ denote $F_p/\Z_2$. We obtain a filtration
$F'_0\subset F'_1\subset \dots \subset F'_{n/2-1}\subset G'_n$ such that the inclusion 
$F'_{n/2-1}\subset G'_n$ is a homotopy
equivalence and 
$$H^j(F'_p, F'_{p-1};\Z_2) \simeq H^{j-2p(m-1)}(\RP^{m-1};\Z_2), \quad p=0, 1, \dots, n/2-1$$
(using the Thom isomorphism). Hence $H^j(F'_p, F'_{p-1};\Z_2)$ is nonzero (and one-dimensional) only 
for $2p(m-1)\le j\le (2p+1)(m-1)$.
The spectral sequence of filtration $F'_p$ 
$$E_1^{p,q} = H^{p+q}(F'_p, F'_{p-1};\Z_2) \Rightarrow H^{p+q}(G'_n;\Z_2)$$
has $E_1^{p,q} \simeq \Z_2$ for $p(2m-3)\le q\le p(2m-3)+(m-1)$ and $E_1^{p,q} =0$ otherwise.
Hence for any differential $d_r$, where $r\ge 1$, either the source or the target vanish.
Therefore, $E_1=E_\infty$ and our statement follows.
\end{proof}

We will now calculate the Stiefel-Whitney classes of the negative normal bundle $\nu_-(V'_p)$.
In particular, we will find out for which $p$ this bundle is orientable. 
This information is needed for computing the integral cohomology of $G_n'$. 

\begin{lemma}\label{lemma2}
 The total Stiefel - Whitney class of the negative normal bundle $\nu_-(V'_p)$
equals
\[(1+\alpha)^{p(m-1)}\, \in \, H^\ast(V'_p;\Z_2),\]
where $\alpha\in H^1(V'_p;\Z_2) \simeq \Z_2$ denotes the generator. 
\end{lemma}
\begin{proof} 
As in the proof of Lemma \ref{lemma1} we obtain that the negative normal bundle $\nu_-(V'_p)$ 
splits as a direct sum
of $2p$ vector bundles $\eta_s$ of rank $m-1$, one for each negative eigenvalue 
$$\lambda_s =2\left[\cos\left(\frac{\pi s}{n+1}\right)-\cos\left(\frac{2\pi k}{n+1}\right)\right]$$ 
of the Hessian. Here $k=n/2 -p$. 

Let $\tau$ denote the tangent bundle of $\RP^{m-1}$. Let $\gamma^\perp$ be a rank $m-1$ vector bundle
over $\RP^{m-1}$ such that its fiber over a line $\ell\in \RP^{m-1}$ is the orthogonal complement $\ell^\perp$.

We claim that 
$$
\eta_s\simeq \left\{
\begin{array}{l}
\tau\quad\mbox{ if $s$ is even,}\\
\gamma^\perp\quad\mbox{ if $s$ is odd.}
\end{array}
\right.
$$
Indeed,
this bundle is obtained from the tangent bundle $\xi$ of $S^{m-1}$ (cf. Lemma \ref{lemma1}) by identifying the
antipodal points, and under this identification the first vector $v_1$ should be replaced by the last vector $v_n$ (cf. (\ref{first})).
Formulae
(\ref{first}) show that 
\[v_n = -\cos(\pi s)\cdot v_1 = (-1)^{s+1} \cdot v_1\]
and hence the bundle $\eta_s$ is obtained from $\xi$ by identifying the fibers over points $a$ and $-a$
with a twist $(-1)^{s+1}$. This implies our claim, cf. \cite{MS}.

For given $p$ there are equal number of negative eigenvalues $\lambda_s$ of the Hessian on $V_p$
with even and with odd $s$.
Therefore the bundle $\nu_-(V'_p)$ is isomorphic to a direct sum of $p$ copies of $\tau \oplus \gamma^\perp$.

The total Stiefel - Whitney class of $\gamma^\perp$ is $(1+\alpha)^{-1}$ and the total Stiefel - Whitney class
of $\tau$ is $(1+\alpha)^m$ (cf. \cite{MS}). Hence the total Stiefel - Whitney class of the negative bundle is 
\[\left[(1+\alpha)^{-1}\cdot (1+\alpha)^m\right]^p = (1+\alpha)^{(m-1)p}.\]
\end{proof}
\begin{corollary}\label{orient1}
If $m$ is odd, then the negative normal bundle $\nu_-(V'_p)$ is orientable for any $p$.

If $m$ is even, then the negative normal bundle $\nu_-(V'_p)$ is orientable for all even $p$ and it is non-orientable 
for all odd $p$.
\end{corollary}
\begin{proof} By Lemma \ref{lemma2} the first Stiefel - Whitney class of  $\nu_-(V'_p)$
is $p(m-1)\alpha$. This implies our statement.
\end{proof}

Recall our permanent assumption $m>1$ and $n$ is even.

\begin{theorem}\label{thm3}
If $m>1$ is odd then
$$
H^j(G'_n;\Z) \simeq \left\{
\begin{array}{ll}
\Z, & \mbox{for $j=2i(m-1)$,\quad where $i=0, 1, \dots, n/2 -1$},\\
\Z_2, & \mbox{for $j$ even satisfying $2i(m-1)<j\le (2i+1)(m-1)$}\\ 
&\mbox{with $i$ as above},\\
0 &\mbox{otherwise}.
\end{array}
\right.
$$
If $m$ is even then
$$
H^j(G'_n;\Z) \simeq \left\{
\begin{array}{ll}
\Z, & \mbox{for $j=(4r+\epsilon)(m-1)$, $r=0, 1, \dots, [(n-2)/4]$, $\epsilon =0, 1$},\\
\Z_2, & \mbox{for $j=4r(m-1)+i$, or $j=(4r'+2)(m-1)+i'$, where}\\
&\mbox{$i=2, 4, \dots, m-2$, $r$ is as above, $i'=1, 3, \dots, m-1$ and}\\
& \mbox{$0\le r'\le (n-4)/4$,}\\
0, &\mbox{otherwise}.
\end{array}
\right.
$$
\end{theorem}

\begin{proof}
Consider filtration $F'_0\subset F'_1\subset \dots \subset F'_{n/2-1}\subset G'_n$ (cf. proof of Corollary \ref{poincare})
and the associated spectral sequence
$${E'}_1^{p,q} = H^{p+q}(F'_p, F'_{p-1};\Z) \Rightarrow H^{p+q}(G'_n;\Z).$$
$F'_p-F'_{p-1}$ contains a single critical submanifold $V'_p\simeq \RP^{m-1}$ with index $2p(m-1)$. The normal bundle to
$V'_p$ is orientable if $p(m-1)$ is even and it is non-orientable if $p(m-1)$ is odd. The Thom isomorphism gives 
\begin{eqnarray}
\quad\quad H^j(F'_p, F'_{p-1};\Z) \simeq \left\{
\begin{array}{ll}
H^{j-2p(m-1)}(\RP^{m-1};\Z), & \mbox{if $p(m-1)$ is even,}\\ \\
H^{j-2p(m-1)}(\RP^{m-1};\pm \Z) & \mbox{if $p(m-1)$ is odd.}
\end{array}
\right.
\end{eqnarray}
Here $\pm Z$ denotes the nontrivial local system of groups $\Z$ over $\RP^{m-1}$; its monodromy
along the generator of $\pi_1(\RP^{m-1})$ is multiplication by $-1$.

For $m$ even we have 
\[H^j(\RP^{m-1};\Z) = \left\{
\begin{array}{ll}
\Z & \mbox{for $j=0$ and $j=m-1$,}\\
\Z_2 & \mbox{for $j=2, \, 4,\, \dots, m-2$},\\
0 & \mbox{otherwise}
\end{array}
\right.
\]
and
\[H^j(\RP^{m-1};\pm \Z) = \left\{
\begin{array}{ll}
\Z_2 & \mbox{for $j=1, \, 3,\, \dots, m-1$},\\
0 & \mbox{otherwise}
\end{array}
\right.
\]
For $m$ odd we have
\[H^j(\RP^{m-1};\Z) = \left\{
\begin{array}{ll}
\Z & \mbox{for $j=0$,}\\
\Z_2 & \mbox{for $j=2, \, 4,\, \dots, m-1$},\\
0 & \mbox{otherwise}.
\end{array}\right.
\]

Therefore, in the above spectral sequence holds
${E'}_1^{p,q} =0$ for $p(2m-3)\le q\le p(2m-3)+(m-1)$.
It implies that for $m>2$, either the source or the target of any differential $d_r$ vanish.

Hence for $m>2$ holds ${E'}_1={E'}_\infty$ and any diagonal $p+q=\const$ contains at most one nonzero
group. This proves our statement for $m>2$.

Assume now that $m=2$ and consider the first differential $d_1:{E'}_1^{r-1,r}\to {E'}_1^{r,r}$.
We have
\[{E'}_1^{r-1,r}\simeq H^{2r-1}(F_{r-1},F_{r-2};\Z) = \left\{
\begin{array}{ll}
\Z, & \mbox{if $r$ is odd,}\\
\Z_2, &\mbox{if $r$ is even}
\end{array}
\right.\]
and 
\[{E'}_1^{r,r}\simeq H^{2r}(F_{r},F_{r-1};\Z) = \left\{
\begin{array}{ll}
0, & \mbox{if $r$ is odd,}\\
\Z, &\mbox{if $r$ is even.}
\end{array}
\right.\]
We see that $d_1$ vanishes since there are no nonzero homomorphisms ${E'}_1^{r-1,r}\to {E'}_1^{r,r}$ for any $r$.

The higher differentials $d_r, \, r\ge 2$ vanish by obvious reasons. Hence, the conclusion we made for $m>2$,
holds also for $m=2$. 
\end{proof}

The following Theorem is the main result of this section.
It describes the multiplicative structure of $H^\ast(G'_n;\Z)$. 
Recall that we assume that $n$ is even.

\begin{theorem}\label{thm4}
For $m>1$ odd, $H^\ast(G'_n;\Z)$ is the commutative ring given by sequence of generators
\[\delta_i\in H^{2i(m-1)}(G'_n;\Z)\simeq \Z, \quad i=0, 1, 2, \dots, \]
and
\[ e\in H^2(G'_n;\Z)\simeq \Z_2,\]
satisfying the following relations
\[
\begin{array}{lll}
\delta_i\delta_j =\displaystyle{\frac{(2i+2j)!}{(2i)!(2j)!}}\cdot \delta_{i+j},\,\,& 
\delta_{n/2}=0,& \delta_0=1,\\ \\
2e=0,& e^{(m+1)/2}=0.&
\end{array}
\]
If $m$ is even then $H^\ast(G'_n;\Z)$ is the graded-commutative ring given by the generators
\begin{eqnarray*}
\delta_i\in H^{4i(m-1)}(G'_n;\Z)\simeq \Z, \quad i=0, 1, 2, \dots, 
\end{eqnarray*}
and also
\[\begin{array}{l}
 e\in H^2(G'_n;\Z)\simeq \Z_2,\quad 
a\in H^{m-1}(G'_n;\Z) \simeq \Z, \quad b\in H^{2m-1}(G'_n;\Z)\simeq \Z_2,
\end{array}
\]
satisfying the following relations
\[
\begin{array}{lll}
\delta_i\delta_j =\displaystyle{\frac{(2i+2j)!}{(2i)!(2j)!}}\cdot \delta_{i+j},\,\, & 
\delta_{[(n+2)/4]}=0,& \delta_0=1,\\ \\
2e=0,& e^{m/2}=0,\\ \\
a^2=0, & ab=0,& ae=0, \\ \\
 2b=0,& b^2=0, \\ \\
\delta_{k}b =0\quad \mbox{(if $n=4k+2$).}&
\end{array}
\]
\end{theorem}
\textit{Remark.} For $m=2$ the generator $e$ disappears since one of the above relations reads $e=0$.
If $m$ is even and $n=2$ then $b=0$ since one of the relations gives $\delta_0 b=0$.

\begin{proof} Consider the universal $\Z_2$-bundle $S^\infty \to \RP^\infty$ and 
the associated fibration $S^\infty\times_{\Z_2} G_n \to
 \RP^\infty$, having $G_n$ as the fiber. The total space $S^\infty\times_{\Z_2} G_n$ is homotopy equivalent to $G_n'$. 
The Serre spectral sequence of this fibration converges to the cohomology algebra $H^\ast(G_n';\Z)$. 
The initial term is
$$E_2^{p,q} = H^p(\RP^\infty;\mathcal H^q(G_n;\Z)),$$
where $\mathcal H^q(G_n;\Z)$, the cohomology of the fiber, is understood as a local system over $\RP^\infty$.

From Theorem \ref{openstring5} we know that $H^q(G_n;\Z)$ is either $\Z$ or trivial.
There are two types of local systems with fiber $\Z$ over $\RP^{\infty}$, which we will denote $\Z$ and $\pm \Z$.
Their structure is determined by the monodromy along any
noncontractible loop of $\RP^\infty$, which is $1$ in the case of $\Z$ and $-1$ in the case of $\pm \Z$. 
 
\begin{figure}[h]
\setlength{\unitlength}{0.9cm}
\begin{center}
\begin{picture}(12,12)
\linethickness{0.3mm}
\put(0,1){\vector(1,0){12}}
\put(0,1){\vector(0,1){11}}

\multiput(0,1)(0,4){3}{\circle*{0.2}}

\multiput(2,1)(2,0){5}{\circle*{0.1}}
\multiput(2,5)(2,0){5}{\circle*{0.1}}
\multiput(2,9)(2,0){5}{\circle*{0.1}}
\multiput(1,3)(2,0){5}{\circle*{0.1}}
\multiput(1,7)(2,0){5}{\circle*{0.1}}
\multiput(1,11)(2,0){5}{\circle*{0.1}}

\multiput(1.2,2.8)(2,0){4}{\vector(3,-2){2.5}}
\multiput(1.2,6.8)(2,0){4}{\vector(3,-2){2.5}}
\multiput(1.2,10.8)(2,0){4}{\vector(3,-2){2.5}}

\put(-1.5,3){$m-1$}
\put(-2,5){$2(m-1)$}
\put(-2,7){$3(m-1)$}
\put(-2,9){$4(m-1)$}
\put(-2,11){$5(m-1)$}

\put(1,0){$1$}
\put(2,0){$2$}
\put(3.7,0){$m+1$}
\put(5.7,0){$m+3$}
\put(7.7,0){$m+5$}
\put(12,0){$p$}
\put(-1,12){$q$}

\end{picture}
\end{center}
\caption{Term $E_m$ of the spectral sequence for $m$ odd}\label{oddd}
\end{figure}
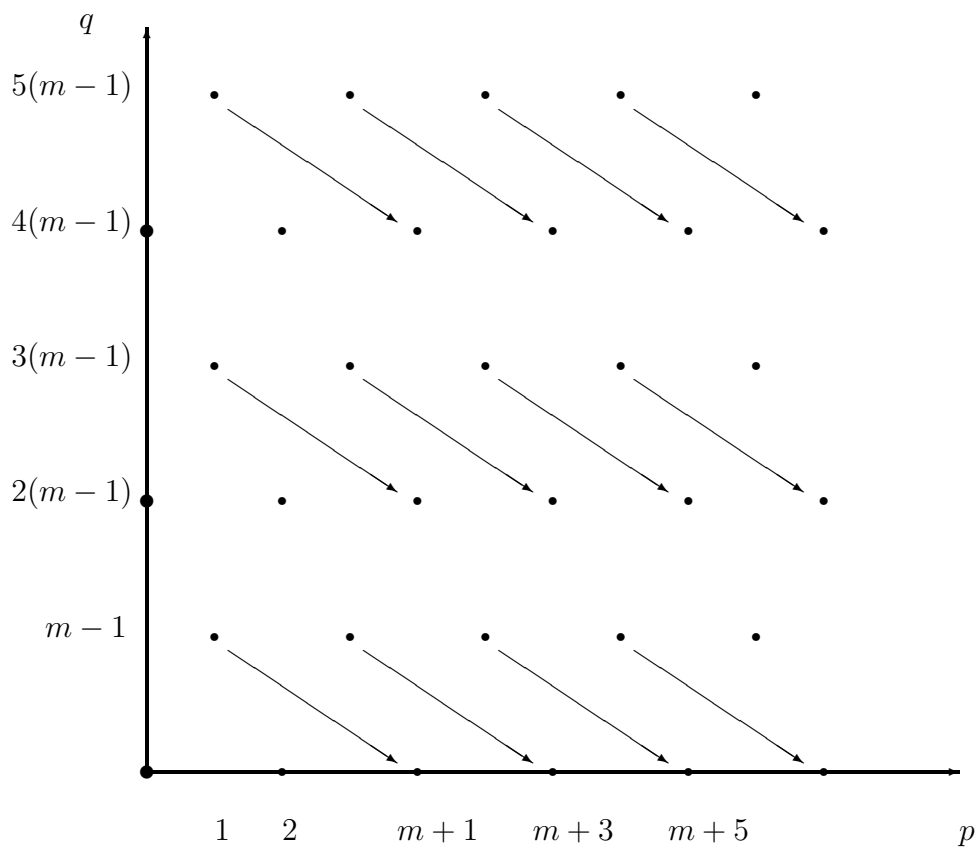

Assume first that $m>1$ is odd. From formula (\ref{reflection2}) we find that 
$$
\mathcal H^q(G_n;\Z) \simeq \left\{
\begin{array}{ll}
\Z, &\mbox{for $q=2i(m-1)$,}\\
\pm \Z &\mbox{for $q=(2i+1)(m-1)$,}
\end{array}
\right.
$$
where $i = 0, 1, \dots, n/2-1$. Hence we find 
$$
E_2^{p,q} = \left\{
\begin{array}{ll}
\Z, &\mbox{for $p=0$ and $q=2i(m-1)$},\\
\Z_2, &\mbox{if $p>0$ is even and $q =2i(m-1)$}\\
&\mbox{or if $p$ is odd and $q=(2i+1)(m-1)$,}\\
0, &\mbox{otherwise},
\end{array}
\right.
$$
where $i = 0, 1, \dots, n/2-1$. As a bigraded algebra, $E_2$ can be identified with the tensor product
$$E_2^{0,\ast} \otimes E_2^{\ast, 0}\otimes A,$$
where 
\begin{eqnarray*}
E_2^{0, \ast} \simeq H^{2\ast}(G_n;\Z),\quad
E_2^{\ast,0} \simeq H^\ast(\RP^\infty;\Z),
\end{eqnarray*}
and $A$ is an exterior algebra with $A^{0,0}\simeq \Z$, and $A^{1,m-1}\simeq \Z_2$. If $x\in E_2^{1, m-1}$ is the 
generator, then relation $x^2=0$ follows from relation 
$\sigma_1^2 =2\sigma_2$ (in the notations of Theorem \ref{openstring5}).
Here we denote $H^{2\ast}(G_n;\Z) \subset H^{\ast}(G_n;\Z)$ the graded subring
$$H^{2\ast}(G_n;\Z)\, =\, \bigoplus_i H^{2i(m-1)}(G_n;\Z).$$
The structure of the ring $H^{2\ast}(G_n;\Z)$
follows from Theorem \ref{openstring5}.

The first nontrivial differential is $d_m$. Since we know the additive structure of $H^\ast(G_n';\Z)$ (cf. Theorem \ref{thm3})
we find that the differential $d=d_m: E_2^{1, m-1}\to E_2^{m+1,0}$ must be an isomorphism. On the other hand 
$d: E_2^{0, 2i(m-1)}\to E_2^{m, (2i-1)(m-1)}$ vanishes (since the range is the zero group).
It follows that $d: E_2^{p, j(m-1)}\to E_2^{p+m,(j-1)(m-1)}$ is nonzero iff both $p$ and $j$ are odd. 

Figure \ref{oddd} shows the nontrivial differential $d=d_m$. 
The fat circles denote group $\Z$ and the small circles denote $\Z_2$.

We conclude that the bigraded  algebra $E_{m+1}$ is isomorphic to the tensor product of algebras
$$H^{2\ast}(G_n;\Z)\otimes H^\ast(\RP^{m-1};\Z),$$ 
where $H^{2i(m-1)}(G_n;\Z)$ has bidegree $(0, 2i(m-1))$ and $H^{2j}(\RP^{m-1};\Z)$ has bidegree $(2j,0)$. 
It is clear that all further differentials vanish and hence $E_\infty = E_{m+1}$. Any diagonal $p+q=c$ contains at most one
nonzero group, and hence the algebra $H^\ast(G'_n;\Z)$ coincides with $E_\infty$.
This proves our statement for $m>1$ odd.

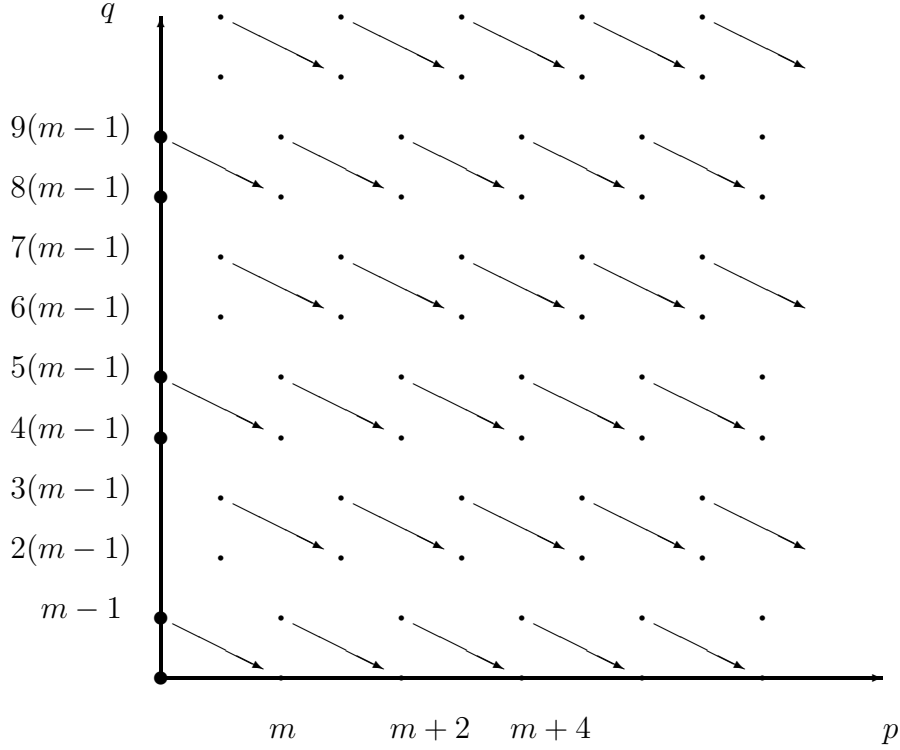
\begin{figure}[h]
\setlength{\unitlength}{0.8cm}
\begin{center}
\begin{picture}(12,12.5)
\linethickness{0.3mm}
\put(0,1){\vector(1,0){12}}
\put(0,1){\vector(0,1){11}}

\multiput(0,1)(0,4){3}{\circle*{0.2}}
\multiput(0,2)(0,4){3}{\circle*{0.2}}

\multiput(2,1)(2,0){5}{\circle*{0.1}}
\multiput(2,2)(2,0){5}{\circle*{0.1}}
\multiput(2,6)(2,0){5}{\circle*{0.1}}
\multiput(2,10)(2,0){5}{\circle*{0.1}}

\multiput(1,4)(2,0){5}{\circle*{0.1}}
\multiput(1,8)(2,0){5}{\circle*{0.1}}
\multiput(1,12)(2,0){5}{\circle*{0.1}}

\multiput(2,5)(2,0){5}{\circle*{0.1}}
\multiput(2,9)(2,0){5}{\circle*{0.1}}

\multiput(1,3)(2,0){5}{\circle*{0.1}}
\multiput(1,7)(2,0){5}{\circle*{0.1}}
\multiput(1,11)(2,0){5}{\circle*{0.1}}

\multiput(0.2,1.9)(2,0){5}{\vector(2,-1){1.5}}
\multiput(0.2,5.9)(2,0){5}{\vector(2,-1){1.5}}
\multiput(0.2,9.9)(2,0){5}{\vector(2,-1){1.5}}

\multiput(1.2,3.9)(2,0){5}{\vector(2,-1){1.5}}
\multiput(1.2,7.9)(2,0){5}{\vector(2,-1){1.5}}
\multiput(1.2,11.9)(2,0){5}{\vector(2,-1){1.5}}

\put(-2,2){$m-1$}
\put(-2.5,3){$2(m-1)$}
\put(-2.5,4){$3(m-1)$}
\put(-2.5,5){$4(m-1)$}
\put(-2.5,6){$5(m-1)$}
\put(-2.5,7){$6(m-1)$}
\put(-2.5,8){$7(m-1)$}
\put(-2.5,9){$8(m-1)$}
\put(-2.5,10){$9(m-1)$}

\put(1.8,0){$m$}
\put(3.8,0){$m+2$}
\put(5.8,0){$m+4$}
\put(12,0){$p$}
\put(-1,12){$q$}
\end{picture}
\end{center}
\caption{Term $E_m$ of the spectral sequence for $m$ even}\label{even}
\end{figure}

Assume now that $m$ is even. Recall that we always assume that $n$ is even.
From formula (\ref{reflection2}) we find that 
$$
\mathcal H^q(G_n;\Z) \simeq \left\{
\begin{array}{ll}
\Z, &\mbox{for $q=4i(m-1)$, or $q=(4i+1)(m-1)$,}\\
\pm \Z, &\mbox{for $q=(4i+2)(m-1)$ or $q=(4i+3)(m-1)$,}
\end{array}
\right.
$$
assuming that $q<n(m-1)$. Hence we find 
$$
E_2^{p,q} = \left\{
\begin{array}{ll}
\Z, &\mbox{for $p=0$ and $q=(4i+\epsilon)(m-1)$, where $\epsilon =0, 1$},\\
\Z_2, &\mbox{if $p>0$ is even and $q =(4i+\epsilon)(m-1)$}\\
&\mbox{or if $p$ is odd and $q=(4i+2+\epsilon)(m-1)$,}\\
0, &\mbox{otherwise}.
\end{array}
\right.
$$
As a bigraded algebra, $E_2$ can be identified with the tensor product
$$E_2^{0,\ast} \otimes E_2^{\ast, 0}\otimes B^{\ast,\ast},$$
where 
\begin{eqnarray*}
E_2^{0, \ast} \simeq H^{4\ast}(G_n;\Z)\otimes C^\ast,\quad
E_2^{\ast,0} \simeq H^\ast(\RP^\infty;\Z),
\end{eqnarray*}
$C^\ast$ is an exterior algebra with $C^0\simeq \Z$ and $C^{m-1} \simeq \Z$
and $B^{\ast,\ast}$ is an exterior bigraded algebra with $B^{0,0}\simeq \Z$, and $B^{1,2(m-1)}\simeq \Z_2$. If $y\in E_2^{1,2(m-1)}$
denotes the generator then $y^2=0$ follows from relation $\sigma_2^2 =2\sigma_4$, cf. Theorem \ref{openstring5}.
We denote by $H^{4\ast}(G_n;\Z) \subset H^{\ast}(G_n;\Z)$ the graded subring
$$H^{4\ast}(G_n;\Z)\, =\, \bigoplus_i H^{4i(m-1)}(G_n;\Z).$$

Consider now the first nontrivial differential $d=d_m:E_2^{p, q}\to E_2^{p+m,q-m+1}$. It is clear that it may be nonzero only
for $q$ of the form $q=(2i+1)(m-1)$.
On the other hand, since we know the additive structure of the limit (cf. Theorem \ref{thm3}), 
we conclude that $d: E_2^{0,m-1}=\Z\to E^{m,0}_2=\Z_2$ is onto. Using multiplicative properties of the spectral sequence
we find that all the differentials shown on Figure \ref{even} are epimorphic. In fact all differentials on Figure \ref{even}, 
except those which start at the $q$ axis, are isomorphisms (since they act between isomorphic groups). As before, the fat circles
denote $\Z$ and small circles denote $\Z_2$.

Hence, to the next term $E_{m+1}$, survive classes $\sigma_{4i}=\delta_i$, and also $a=2\sigma_1$,
$e\in E_{m+1}^{2,0}\simeq \Z_2$ and $b\in E_{m+1}^{1,2(m-1)}$ and their products $\delta_ia$, $\delta_ie^j$,
 and $\delta_ibe^j$ with
$j<m/2$. It is clear that all further differentials vanish and in each diagonal $p+q=c$ there is at most one nonzero group.
Therefore we conclude that ring $H^\ast(G'_n;\Z)$ is isomorphic to $E_{m+1}$. Its structure coincides with the description given in Theorem \ref{thm4}.
\end{proof}

\section{\bf   Equivariant Lusternik - Schnirelman theory via nonsmooth critical point theory}

In this section we firstly recall the basic notions of the critical point theory for nonsmooth
functions, suggested recently in \cite{C} and \cite{DM}.
Secondly, we apply the nonsmooth critical point theory to get a simple independent exposition of a version of the equivariant Lusternik - Schnirelman theory of \cite{Mar}, \cite{CP}, which we need for our applications to the billiard problems. 
One of the advantages of our approach is applicability to manifolds with boundary.

Let $X$ be a metric space endowed with the metric $d$. Given a point $p\in X$ and 
$\delta>0$, we denote by $B(p, \delta)\subset X$ the ball of radius $\delta$ centered at $p$.

\begin{definition}\label{slope}
 Let $f: X\to \R$ be a continuous function. 
{\it The weak slope of $f$ at} a point 
$p\in X$, denoted $|df|(p)$, is defined as  
the supremum of all $\sigma \in [0, \infty]$ such that there exists $\delta>0$ and a continuous deformation 
$\eta: B(p, \delta)\times [0, \delta] \to X$ with the following properties
$$
\begin{array}{cc}
d(\eta(q,t),q) \le t, & \quad f(\eta(q,t)) \le f(q) - \sigma t
\end{array}
$$
for all $q\in B(p, \delta)$, $t\in [0,\delta]$. 

A point $p\in X$ is said to be {\it a critical point of function $f$} if $|df|(p)=0$. 
\end{definition}

\begin{example} Let $X$ be a smooth Riemannian manifold without boundary, and let $f: X\to \R$ be a smooth function.
Then the weak slope $|df|(p)$ coincides with the norm of the differential $||df(p)||$, 
viewed as a bounded linear functional on the tangent space $T_p(X)$. 
\end{example}

\begin{example} Let $X$ be a smooth Riemannian manifold with boundary, and let $f: X\to \R$ be a smooth function.
A point on the boundary $p\in \partial X$ is a critical point of $f$ iff there is no tangent vector $v\in T_pX$
pointing inside $X$ such that the derivative $v(f)<0$ is negative. 
The last condition implies that 
\begin{eqnarray}\label{grad}
df_p|_{T_p\partial X}=0,
\end{eqnarray}
i.e. the gradient of $f$ at point $p\in \partial X$ 
is orthogonal to the boundary $\partial X$. A point $p\in \partial X$ is a critical point of $f$ iff 
(\ref{grad}) holds and the gradient of $f$ at $p$ has the inward direction. 
It is clear that the above conditions are independent of the Riemannian metric.
\end{example}

\begin{proposition}\label{ls} Let $f:X\to \R$ be a continuous function on a compact metric space
$X$. Then the number of critical points of $f$ (in the sense of Definition \ref{slope})
is at least $\cat(X)$, the Lusternik - Schnirelman category of $X$.
\end{proposition}
 
This follows from a much more general theorem 3.7 of  \cite{DM}.

We will apply the nonsmooth critical point theory to equivariant critical point theory of smooth functions.
Compare \cite{Mar}, \cite{CP}. 

\begin{proposition}\label{prop16}
Let $M$ be a smooth compact 
$G$-manifold with boundary, where $G$ is a finite group.
Let $f: M\to \R$ be a  $G$-invariant smooth function. 
Suppose that at points of the boundary $p\in \partial M$ 
the gradient of $f$ does not vanish and has the outward direction. Then the number of $G$-orbits of points
$p\in M$ with $df_p=0$ is at least $\cat(M/G)$.
\end{proposition}

\begin{proof} Let $X$ denote the space of orbits $X=M/G$. Function $f$
determines a continuous function $\tilde f: X\to \R$.  We want to show that any orbit $x\in X$, representing
points $p\in X$ with $df_p\ne 0$, is not a critical point of $\tilde f: X\to \R$ 
in the sense of Definition \ref{slope}. This would imply that the number of critical orbits of $f$ is at least the
number of critical points of $\tilde f$; the later can be estimated from below by $\cat(X)$ 
by Proposition \ref{ls}.

We will assume that $M$ is supplied with a $G$-invariant Riemannian metric.
Let $p\in M$ be a point with $df_p\ne 0$. 
We want to construct a smooth vector field $v$ in a neighborhood
of the orbit of $p$ having the following properties:
\begin{itemize}
\item[(a)] $v(f)_p <0$;
\item[(b)] the norm of vector $v_p$ equals 1;
\item[(c)] $v$ is $G$-invariant;
\item [(d)] in case $p$ belongs the boundary $\partial M$, the vector $v_p$ points inside $M$.
\end{itemize}
To construct such vector field $v$ one first finds a vector $v_p\in T_pM$ for each point $p$
of the orbit so that (a), (b) and (d) are satisfied. 
It is then possible to extend the vectors $v_p$ to form a smooth
vector field $v$ in a neighborhood of the orbit of $p$ with properties (a), (b), (d). Then (c) can be achieved 
by averaging. 

The flow determined by the vector field $v$, gives a deformation of a ball around the point $\{p\}\in X$, which represents the orbit of $p$, showing that
the slope $|d\tilde f|(p)$ is positive. 
\end{proof}

\section{\bf Proof of Theorem \ref{thm2}}

Let $T\subset \R^{m+1}$ be a compact strictly convex domain bounding a smooth hypersurface $X=\partial T$.
Given a point $A\in X$ and an integer $n$, consider the configuration space $G_n=G(X;A,A,n)$ (cf. (\ref{closedstring})) and 
the smooth function
\begin{eqnarray*}
L_X : G_n\to \R, \quad L_X(x_1, \dots, x_n) = -\sum\limits_{i=0}^{i=n}|x_i -x_{i+1}|
\end{eqnarray*}
(the negative total length), where we understand $x_0=A=x_{n+1}$. This function is invariant with respect to
reflection $T: G_n\to G_n$ (cf. (\ref{reflection})). Hence 
$L_X$ determines a continuous function 
\begin{eqnarray}\label{lprime}
L'_X: G'_n\to \R, \quad G'_n=G_n/T,
\end{eqnarray}
and the critical points 
of $L'_X$ (in the sense of Definition \ref{slope}) are in one-to-one correspondence with $\Z_2$-orbits of 
closed billiard trajectories in $X$, which start and end at $A$ and make $n$ reflections. This follows from 
Lemma 2 in \cite{F} and from the argument in the proof of Proposition \ref{prop16}.

Note that for $n$ is even, $T$ acts freely on $G_n$, the 
factor $G'_n=G_n/T$ is a smooth manifold and the function $L'_X$ is smooth. In this case the nonsmooth critical
point theory coincides with the usual one. 

We claim:

{\it The number of critical points of $L'_X$ is at least the Lusternik - Schnirelman category $\cat(G'_n)$;
moreover, assuming that $n$ is even and the function 
$L'_X$ is Morse, the number of critical points of (\ref{lprime}) is at least the sum of Betti
numbers of $G'_n$ with $\Z_2$ coefficients.} 

The italiced statement does not follow directly from the traditional Morse-Lusternik-Schnirelman 
theory, since $G'_n$ is not compact. However, as in \cite{FT}, \cite{F}, one may fix  $\epsilon>0$
and consider the following compact subset 
$$G_\epsilon\subset G_n, \quad G_\epsilon =\{(x_1, \dots, x_n)\in X^{\times n}: \prod\limits_{i=0}^{n} |x_i-x_{i+1}| \ge \epsilon\}.
$$
If $\epsilon>0$ is small enough then: 

\begin{itemize}
\item[(a)] $G_\epsilon$ is a compact manifold with boundary;

\item[(b)]  the inclusion $G_\epsilon\subset G_n$ is a $\Z_2$-equivariant homotopy equivalence;

\item[(c)]  all critical points of $L_X$ are contained in $G_\epsilon$;

\item[(d)]  at every point of $\partial G_\epsilon$ the gradient of $L_X$ has the outward direction.
\end{itemize}
Cf. \cite{FT}, Proposition 4.1 and \cite{F}, Proposition 4.

Since $G'_\epsilon = G_\epsilon/T$ is a compact smooth manifold with boundary, 
we may apply the Morse-Lusternik-Schnirelman theory
to it. Condition (d) implies that the critical points of the restriction of $L'_X$ on $\partial G'_\epsilon$
should not be taken into the account (cf. Proposition \ref{prop16}). 
Therefore, the number of critical points of $L'_X|_{G'_\epsilon}$
is at least $\cat(G'_\epsilon)= \cat(G'_n)$. If $L'_X|_{G'_\epsilon}$ is Morse, then the number of its critical points is at least
the sum of Betti numbers of $G'_\epsilon$, which is the same as the sum of Betti numbers of $G'_n$. 

In the proof of statement (I) of Theorem \ref{thm2} we will use the following general simple remark:

{\it for any regular covering map $p: \tilde X\to X$ holds}
\begin{eqnarray}\label{cover}
\cat(X) \ge \cat(\tilde X).
\end{eqnarray}
Indeed, if $A\subset X$ is an open subset which is contractible to a point in $X$ then 
$\tilde A =p^{-1}(A)$ is a disjoint union of open subsets of $\tilde X$, such that each is contractible to a point in $\tilde X$. 
Hence, any categorical open cover $A_1\cup A_2\cup\dots \cup A_k = X$ produces a categorical open cover
 $\tilde A_1\cup \tilde A_2\cup\dots \cup \tilde A_k$ of $\tilde X$.

The above remark applies to the two-fold cover $G_n \to G'_n$ giving
\[\cat(G'_n) \ge \cat(G_n) \ge \cl(G_n) +1,\]
where $\cl(G_n)$ is the cohomological cup-length of $G_n$ (Froloff - Elsholz Theorem).

Now we use Theorem \ref{openstring5} to compute the cup-length of $G_n$. In $m\ge 3$ is odd then 
$\sigma_1^{n-1}=(n-1)! \sigma_{n-1} \ne 0\in H^{(n-1)(m-1)}(G_n;\Z)$ and hence $\cl(G_n)=n-1$. Therefore
 $\cat(G'_n)\ge n$.
This proves statement (I) in the case $m\ge 3$ odd.

If $m$ is even, then the longest nontrivial cup-product in $H^\ast(G_n;\Z)$ 
is 
$$
\left\{
\begin{array}{ll}
\sigma_1\sigma_2^{k-1}&\mbox{if $n=2k$ is even},\\ 
\sigma_2^k,&\mbox{if $n=2k+1$ is odd.}
\end{array}\right.
$$ 
We conclude that for $m$ even the cup-length of $G_n$ equals $[n/2]$ and therefore $\cat(G'_n)\ge [n/2]+1$. 
Together with the information collected above
this proves statement (I) for $m$ even. 

To prove statement (II) of Theorem \ref{thm2} we will use Theorem \ref{thm4} to estimate the Lusternik - Schnirelman
category of $G'_n$. Note that Theorem \ref{thm4} requires the assumption that $n$ is even. Suppose first that $m\ge 3$ is odd. 
Then (in the notations of Theorem \ref{thm4}) we have a nonzero cohomology product
\begin{eqnarray}\label{prod}
\delta_1\delta_2\delta_{2^2}\dots \delta_{2^s}e^{(m-1)/2},
\end{eqnarray}
where $s$ is the largest integer with $2^{s+1}-1\le n/2 -1$, i.e. $s = [\log_2(n)]-2$. 
Note that class $e\in H^2(G'_n;\Z)$ has order $2$, i.e.
$2e=0$. Nontriviality of the above product is equivalent to the claim 
that the product $\delta_1\delta_2\delta_{2^2}\dots \delta_{2^s}$ is an odd multiple of the class
$\delta_{1+2+2^2+\dots+2^s}$. Indeed, we use the relation
$$
\delta_i\delta_j = \left(\begin{array}{c}
2i+2j\\
2i
\end{array}
\right) \delta_{i+j},
$$
and the well-known fact that the binomial coefficient $\left(\begin{array}{c}
2i+2j\\
2i
\end{array}\right)
$
 is even if and only if $i$ and $j$, in their binary expansions,
have a 1 on the same place.

Now we will use the notion of {\it category weight of a cohomology class}, introduced by E. Fadell and S. Huseini \cite{FH}.
We claim that the category weight of $e\in H^2(G'_n;\Z)$ equals $2$.
Indeed, $e$ equals the image of a class $e'\in H^1(G'_n;\Z_2)$ under the Bockstein homomorphism 
$\beta: H^1(G'_n;\Z_2)\to H^2(G'_n;\Z)$, i.e. $e=\beta(e')$. Hence by Theorem (1.2) of E. Fadell and S. Huseini \cite{FH},
the category weight of $e$ is $2$. Therefore, nontriviality of product (\ref{prod}) implies 
\[\cat(G'_n) \ge (s+1) +2\cdot \frac{m-1}{2} +1 = [\log_2n]+m-1.\]
This proves statement (II) for $m\ge 3$ odd.

Consider now the case when $m\ge 2$ is even. 
First we will assume that $n\ge 8$ and $n+2$ is not a 
power of $2$.
Then we have a nontrivial cohomological product
\begin{eqnarray}\label{prod3}
\delta_1\delta_2\delta_{2^2}\dots \delta_{2^s}be^{(m-2)/2},
\end{eqnarray}
where $s=\left[\log_2 \left[\displaystyle{\frac{n+2}{4}}\right]\right]-1$. As above, nontriviality of (\ref{prod1}) implies 
\[\cat(G'_n) \ge (s+1) +1+2\cdot \frac{m-2}{2} +1 = \left[\log_2 \left[\displaystyle{\frac{n+2}{4}}\right]\right]
+m\ge [\log_2n] +m-2.\]

If $n+2=2^r$ is a power of $2$, where $r\ge 4$, then the product 
\begin{eqnarray*}\label{prod4}
\delta_1\delta_2\delta_{2^2}\dots \delta_{2^s}e^{(m-2)/2},
\end{eqnarray*}
(we skip $b$ because of the last relation in Theorem \ref{thm4}) is nonzero, where $s=r-3$.
In this case we obtain
\[\cat(G'_n) \ge (s+1)+2\cdot \frac{m-2}{2} +1 =s+m= [\log_2n] +m-2.\]

We are left to consider the cases $n= 2, \, 4,\, 6$ with $m$ even. Here we have a nontrivial cup-product $be^{(m-2)/2}$
and hence
\[\cat(G'_n) \ge 1+2\cdot \frac{m-2}{2} +1 =m.\]
This implies the estimate of statement (II) of Theorem \ref{thm2} for the specified values of $n$ and $m$. 

For $m>1$ statement (III) of Theorem \ref{thm2} follows from Corollary \ref{poincare}.
If $m=1$ then the space $G'_1=G(S^1;A,A,n)/\Z_2$ consists of $n/2$ connected components and each is contractible;
this can be established by arguments similar to those used in \S 7 of \cite{F}. Hence the sum of Betti numbers of $G'_1$
is $n/2$.

\qed

\section{\bf Cohomology of cyclic configuration spaces of spheres}

The cyclic configuration space $G(X,n)$ of a space $X$ is defined (cf. \cite{FT}) as the set of all configurations 
$(x_1, \dots, x_n)\in X^{\times n}$ with $x_i\ne x_{i+1}$ for $i=1, \dots, n-1$ and $x_n\ne x_1$.
The dihedral group $D_n$ acts naturally on $G(X,n)$. 

In \cite{FT} we showed that information about the cohomology ring of the factor-space $G(S^m,n)/D_n$
leads to estimates on the number of $n$-priodic orbits of convex billiards in $(m+1)$-dimensional space $\R^{m+1}$.
 The rings $H^\ast(G(S^m,n);\Z_2)$ and  $H^\ast(G(S^m,n)/D_n;\Z_2)$ were computed in \cite{FT}.

In this section we will describe the cohomology of the cyclic configuration space $G(S^m,n)$ 
with other fields of coefficients. It turns out that the answer depends on the parity of $m$; therefore we state
the even and the odd dimensional cases in the form of two separate theorems.

The results of this section will be used in the proof of Theorem \ref{thm5}.

\begin{theorem}\label{openstring6}
Let  $m\ge 3$ be odd.
 The ring $H^\ast(G(S^m,n);\Q)$ is given by generators 
\[u\in H^m(G(S^m,n);\Q), \quad \sigma_i \in \, H^{i(m-1)}(G(S^m,n);\Q), \quad i=1, \dots , n-2,\]
and relations:
\begin{eqnarray}
u^2 =0, \quad \sigma_i\sigma_j \, =\, \left \{
\begin{array}{ll}
\displaystyle{\frac{(i+j)!}{i!\cdot j!}}\cdot \sigma_{i+j},&\mbox{if}\quad i+j\le n-2,\\ \\
0, &\mbox{if}\quad  i+j>n-2, 
\end{array}
\right .
\end{eqnarray}
\end{theorem}

One may show that the statement of Theorem \ref{openstring6} holds with $\Q$ replaced by an 
arbitrary field of coefficients $\kk$. However,
the short proof we give below, works only in the case $\kk=\Q$. On the other hand for our purposes of this paper it is enough to 
know the rational cohomology $H^\ast(G(S^m,n);\Q)$. The case of a field $\kk$ of positive characteristic 
may be proven using Theorem
3 of \cite{FT} and computing the spectral sequence similarly to the proof of Theorem 4 in \cite{FT}.

The following theorem gives the answer for $m$ even.

\begin{theorem}\label{openstring8}
Let $\kk$ be a field of characteristic $\ne 2$.
For any even $m\ge 2$ and odd $n\ge 3$  the cohomology algebra $H^\ast(G(S^m,n);\kk)$ has generators 
\[w\in H^{2m-1}(G(S^m,n);\kk), \quad 
\sigma_{2i} \in \, H^{2i(m-1)}(G(S^m,n);\kk),
\quad  i= 1, \dots , (n-3)/2,\]
such that
\begin{eqnarray}\label{relnew}
 w^2 =0, \quad \sigma_{2i}\sigma_{2j} \, =\, \left \{
\begin{array}{ll}
\displaystyle{\frac{(i+j)!}{i!\cdot j!}}\cdot \sigma_{2(i+j)}, &\mbox{if $i+j\le (n-3)/2$.}\\ \\
0, & \mbox{if} \quad i+j > (n-3)/2.
\end{array}
\right .
\end{eqnarray}
\end{theorem}

\begin{proof}[\bf Proof of Theorem \ref{openstring6}]
Consider fibration
\begin{eqnarray}\label{fibration}
p: G(S^m,n)\to S^m,
\end{eqnarray}
where the image of a cyclic configuration $(x_1, \dots, x_n)\in G(S^m,n)$ under projection $p$
is given by
$p(x_1, \dots, x_n) = x_1$. The fiber of $p$ is the configuration space $G(S^m;A,A,n-1)$.
Consider the Serre's spectral sequence of this fibration.
The cohomology of the fiber $ G(S^m;A,A,n-1)$ is described by Theorem \ref{openstring5}; it has generators
$\sigma_1, \dots, \sigma_{n-2}$, which multiply according to (\ref{prod1}). 

This spectral sequence may
have only one nonzero differential $d_m$. We will show that this differential vanishes $d_m=0$.
This would clearly imply our statement.

Since we may write $\sigma_i = (i!)^{-1}(\sigma_1)^i$, it is enough to
show that $d_m(\sigma_1)=0$. 
Vanishing  $d_m(\sigma_1)=0$ follows from the fact that fibration (\ref{fibration}) admits a continuous section
$s: S^m \to G(S^m,n)$ and thus the transgression is trivial. 
To construct $s$, fix a nowhere zero tangent vector field $V$ on the sphere $S^m$ 
(recall that $m$ is odd).
For $x\in S^m$, the tangent vector $V(x)$ determines a half circle starting at $x$, tangent to $V(x)$
and ending at the antipodal point $-x$. Then the section $s$ can be defined by
\[s(x) = (x_1, x_2, \dots, x_n),\quad x\in S^m,\]
where $x_1=x, x_n = -x$ and the points $x_2, \dots, x_{n-1}$ are situated on the half circle making equal angles
as shown on the following
picture.

\begin{figure}
\centerline{\includegraphics{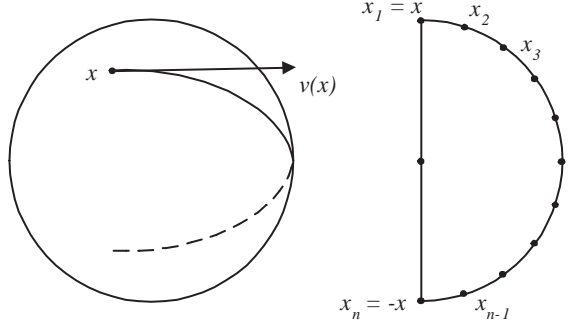}}
\caption{Continuous family of configurations}
\end{figure}
Analytically we may write
\[x_j = \cos\left(\frac{(j-1)\pi}{n-1}\right) x + \sin\left(\frac{(j-1)\pi}{n-1}\right) V(x), \quad j=1, \dots, n.\]
\end{proof}

\begin{proof}[\bf Proof of Theorem \ref{openstring8}]

First we will assume that $m>2$; the case $m=2$ will be treated separately later.

We will describe the additive structure of $H^\ast(G(S^m,n);\kk)$, using approach of the Morse theory.

Let $S^m\subset \R^{m+1}$ be the unit sphere. Consider the total length function
\[L: G(S^m,n) \to \R,\]
where for $(x_1, \dots, x_n)\in G(S^m,n)$ we have
\[L(x_1, x_2, \dots, x_n) = -|x_1-x_2| -|x_2-x_3|-\dots -|x_n-x_1|.\]
The critical points of $L$ are $n$-periodic billiard trajectories in the unit sphere; hence
the critical configurations are 
regular $n$-gons lying in two-dimensional central sections of the sphere. 
A regular $n$-gon is determined by two first vectors $x_1, x_2\in S^m$, which must make an angle
of the form
\[\alpha_p = \frac{2\pi}{n}\cdot (\frac{n-1}{2}-p), \quad\mbox{where}\quad p=0, 1, \dots, (n-3)/2.\]
Recall, that we assume that $n$ is odd.
Fixing $p=0, 1, \dots, (n-3)/2$, we obtain a variety of critical configurations, which we will denote by
 $V_p\subset G(S^m,n)$.
Each $V_p$ has dimension $2m-1$ and is diffeomorphic to the Stiefel manifold 
of pairs of mutually orthogonal vectors in $\R^{m+1}$. Since we assume that $m$ is even and the characteristic
of $\kk$
is  $\ne 2$, we have $H^\ast(V_p;\kk) \simeq H^\ast(S^{2m-1};\kk)$. Note also that $V_p$ is simply connected
(since $m>2$).

I.K. Babenko has shown (cf. \cite{Ba}, Proposition 3.1) 
that function $L$ is nondegenerate in the sense of Bott and the
index of each critical submanifold $V_p$ equals $2p(m-1)$.
Moreover, it is clear that $L(V_p) < L(V_{p'})$ for $p<p'$.

Fix $\varepsilon>0$ small enough and consider the submanifold
$G_{\varepsilon} (S^m,n) \subset G(S^m,n)$, where 
\begin{eqnarray*}
G_{\varepsilon} (S^m,n) = \{(x_1, \dots, x_n) \in (S^m)^{\times n} :
\prod_{i=1}^n|x_i - x_{i+1}| \geq \varepsilon \}.
\end{eqnarray*}
If $\varepsilon >0$ is sufficiently small then (according to Proposition 4.1 of \cite{FT})
$G_{\varepsilon} (S^m,n)$ is a compact manifold with boundary containing all the critical points of $L$
and such that the inclusion $G_{\varepsilon} (S^m,n) \subset G(S^m,n)$ is a $D_n$-equivariant homotopy
equivalence. Moreover, at every point of the boundary $\partial G_{\varepsilon} (S^m,n)$ the gradient
of $L$ has the outward direction.

Choose constants $c_0, c_1, \dots, c_{(n-3)/2}$ such that $L(V_p)<c_p<L(V_{p+1})$ for $0\le p <(n-3)/2$
and $c_{(n-3)/2}=0$. 
Let 
$$F_p=L^{-1}((-\infty, c_p])\cap G_{\varepsilon}(S^m,n).$$
We obtain a filtration 
$F_0\subset F_1\subset \dots\subset F_{(n-3)/2}=  G_\varepsilon (S^m,n).$
Since the inclusion $F_{(n-3)/2}\subset G(S^m,n)$ is a homotopy equivalence, 
we may use the spectral sequence of this filtration in order to calculate the
cohomology of $G(S^m,n)$.

We claim that {\it this filtration is perfect}, i.e. the Poincar\'e polynomial of 
the cyclic configuration space $G(S^m,n)$ equals the sum of the Poincar\'e polynomials of the pairs
$(F_p,F_{p-1})$. The initial term of the spectral sequence is 
\[E^{p,q}_1 = H^{p+q}(F_p,F_{p-1};\kk).\]
Uisng the Thom isomorphism (recall that $V_p$ is simply connected), 
we find that $H^j(F_p,F_{p-1};\kk)$ is isomorphic to 
$H^{j-2p(m-1)}(V_p;\kk)$; hence $H^j(F_p,F_{p-1};\kk)$ is one-dimensional for $j=2p(m-1)$ and for
$j=2p(m-1)+2m-1$ and vanishes for all other values of $j$. This follows since in $F_{p+1}-F_p$ there is
a single non-degenerate
critical submanifold $V_p$, which has index $2p(m-1)$ and $V_p$ is diffeomorphic to the Stiefel manifold 
$V_{m+1, 2}$.

The gradient of $L$ at points of the boundary $\partial G_{\varepsilon}(S^m,n)$ has the outward direction and hence the
points of the boundary do not contribute to the usual statements of the Morse - Bott critical point theory.

For a given $p$ there are precisely two values of $q$ such that $E_1^{p,q}$ is nonzero
($q=p(2m-3)$ and 
$q=p(2m-3)+2m-1$).
From the geometry of the differentials we see that all the differentials $d_r$, $r\ge 1$, 
must vanish if $m>2$. This proves that the cohomology 
$H^j(G(S^m,n);\kk)$ is one dimensional for 
$j=2p(m-1)$ and $j=2p(m-1)+2m-1$, where $p=0, 1, \dots, (n-3)/2,$ and $H^j(G(S^m,n);\kk)$
vanishes for other values of $j$.

Having recovered the additive structure of $H^\ast(G(S^m,n);\kk)$, we may use Theorem \ref{openstring5} to find its 
multiplicative structure. The mapping $(x_1, x_2, \dots, x_n)\mapsto x_1$
is a Serre fibration
$G(S^m,n)\to S^m$; its fiber is $G(S^m:A,A,n-1)$. The Serre spectral sequence has only two nonzero columns
and $d_m$ is the only differential which could be nonzero. 
In the $0$-th columns we have classes $\sigma_i$ in dimensions $i(m-1)$, and in the 
$m$-th column we have classes $\sigma_i u$ having dimension $i(m-1)+m$, cf. Theorem \ref{openstring5}.
Since we already know the additive structure of $H^\ast(G(S^m,n);\kk)$, we conclude that
the differential 
\[d_m: E^{0,i(m-1)}_2\to E_2^{m,(i-1)(m-1)}\]
is an isomorphism for $i$ odd and vanishes for $i$ even. Hence the classes 
$\sigma_{2i}$ and $\sigma_{2i+1}u$ survive. Now, we set $w=\sigma_1u$,
and the conclude that $H^\ast(G(S^m,n);\kk)$ has the multiplicative structure as stated in Theorem \ref{openstring8}.

\setlength{\unitlength}{0.7cm}
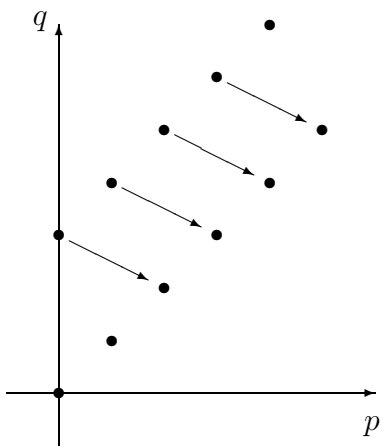
\begin{figure}[h]
 \begin{center}
\begin{picture}(8, 9)
\linethickness{0.15mm}
\put(0,1){\vector(1,0){7}}
\put(1,0){\vector(0,1){8}}
\put(1,1){\circle*{0.2}}
\put(1,4){\circle*{0.2}}
\put(2,5){\circle*{0.2}}
\put(2,2){\circle*{0.2}}
\put(3,3){\circle*{0.2}}
\put(3,6){\circle*{0.2}}
\put(4,4){\circle*{0.2}}
\put(4,7){\circle*{0.2}}
\put(5,5){\circle*{0.2}}
\put(5,8){\circle*{0.2}}
\put(6,6){\circle*{0.2}}
\put(6.8,0.3){$p$}
\put(0.5,8){$q$}
\linethickness{0.05mm}
\put(1.2,3.9){\vector(2,-1){1.5}}
\put(2.2,4.9){\vector(2,-1){1.5}}
\put(3.2,5.9){\vector(2,-1){1.5}}
\put(4.2,6.9){\vector(2,-1){1.5}}
\end{picture}
\end{center}
\caption{Nonzero terms of the spectral sequence for $m=2$}
\end{figure}

For $m=2$ the above argument, based on the spectral sequence of the filtration 
$F_0\subset F_1\subset \dots\subset F_{(n-3)/2}=  G_\varepsilon (S^m,n),$
is not sufficient, since, in principle, this spectral sequence could have a nonzero differential,
shown on the picture. Also, for $m=2$ the critical submanifolds $V_p$ are not simply connected
and so the Thom isomorphisms for the negative normal bundles of the Hessian may require additional twists
by flat line bundles (depending on the orientability of the negative normal bundles of the critical submanifolds).

However, in the case $m=2$ a different argument can be applied.
Consider the action of $SO(3)$ on $G(S^2,n)$ arising from the standard
action of $SO(3)$ on $S^2$. 
Fix a point $A\in S^2$ and consider $G(S^2;A,A,n-1)$ as being canonically embedded in $G(S^2,n)$. We
obtain the map 
\begin{eqnarray}\label{fibration1}
SO(3)\times G(S^2;A,A,n-1) \to G(S^2,n),\quad (R,c)\mapsto Rc,
\end{eqnarray}
given by applying an orthogonal matrix $R\in SO(3)$ to a configuration of points on the sphere 
$c\in G(S^2;A,A,n-1)$.
It is easy to see that (\ref{fibration1}) is a fibration with 
fiber $S^1$. If $c=(A,x_1,\dots, x_{n-1})$ is a configuration of points on $S^2$ such that $A\ne x_1$, 
$x_i\ne x_{i+1}$ for $i=1, \dots, n-1$ and $x_{n-1}\ne A$, then the fiber of fibration (\ref{fibration1}) over $c$ consists
of the space of all pairs $(R_{-\phi}, R_{\phi}(c))$, where $R_{\phi}\in SO(3)$ denotes the rotation by angle $\phi\in [0,2\pi]$
about $A$.

The cohomology algebra of the total space of this fibration
$$H^\ast(SO(3)\times G(S^2;A,A,n-1);\kk)\simeq H^\ast(S^3;\kk)\otimes H^\ast(G(S^2;A,A,n-1);\kk)$$
is given by Theorem \ref{openstring5}. It
has a generator $w$, with $\deg w=3$
(coming from a generator of $H^3(SO(3);\kk)$) and also classes
$\sigma_i$, where $i=0, 1, \dots, n-2$, with $\deg \sigma_i =i$, which are pullbacks of the generators
of $H^\ast(G(S^2;A,A,n-1);\kk)$, cf. Theorem \ref{openstring5}. We have the relation
$w^2=0$ and each product $\sigma_i\sigma_j$
equals a multiple of $\sigma_{i+j}$, the coefficient indicated in formula (\ref{prod2}).

Let us show that {\it the restriction map from the total space to the fiber
\[H^1(SO(3)\times G(S^2;A,A,n-1);\kk)\to H^1(S^1;\kk)\]
is onto}. Since $H^1(SO(3);\kk)=0$, our statement is equivalent to the following. 
Let $c=(A, x_1, \dots, x_{n-1})$ be a fixed configuration. We obtain an embedding
$f: S^1\to G(S^2;A,A,n-1)$ given by $\phi\mapsto R_{\phi}(c)$, where $\phi \in [0,2\pi]$. 
We claim that the induced map 
$f^\ast: H^1(G(S^2;A,A,n-1);\kk) \to H^1(S^1;\kk)$ is onto. In other words, we want to show that the cohomology class
$f^\ast(\sigma_1)\in H^1(S^1;\kk)$ is nonzero.

We may assume that the antipode $A'$ of $A$
does not appear in the configuration $c$. Identify $S^2-A'$ with $\R^2$ using the stereographic projection
with $A'$ as a center; this leads to the following commutative diagram
$$
\begin{array}{clc}
S^1 & \stackrel g\longrightarrow & G(\R^2;0,0,n-1)\\ \\
& \searrow\, \,  f &\downarrow h\\ \\
& & G(S^2;A,A,n-1)
\end{array}
$$
where $g$ is given by rotations of a fixed configuration $c'=(0, y_1, y_2, \dots, y_{n-1})$ of points on the plane,
$c'\in G(\R^2;0,0,n-1)$,
around the origin $0\in \R^2$.
Clearly, the space $G(\R^2;0,0,n-1)$ is homotopy equivalent to $G(\R^2,n)$ and thus the cohomology algebra
$H^\ast(G(\R^2;0,0,n-1);\kk)$, as given by Proposition 2.2 of \cite{FT}, has 1-dimensional generators 
$s_1, \dots, s_n$, which satisfy the relations $s_i^2=0$ for $i=1, \dots, n$ and also a relation of degree
$n-1$, cf. formula (4) in \cite{FT}. From {\it Remark 9} in \cite{F}
we obtain that
\begin{eqnarray}\label{restric}
h^\ast(\sigma_1) = \sum_{i=1}^n (-1)^{i+1}s_i.
\end{eqnarray}

Let $s\in H^1(S^1;\kk)$ denote the generator corresponding to the usual
anti-clockwise orientation of the circle. Then
\begin{eqnarray}\label{restric1}
g^\ast(s_i) = s, \quad i=1, 2, \dots, n.
\end{eqnarray}
Indeed, $g^\ast(s_i) = d_i s$, where $d_i$ is the degree of the following map $S^1\to S^1$
\[\phi\mapsto \frac{R_\phi(y_i) -R_\phi(y_{i-1})}{|R_\phi(y_i) -R_\phi(y_{i-1})|} = 
R_\phi(\frac{y_i -y_{i-1}}{|y_i -y_{i-1}|}), \quad \phi\in [0, 2\pi],\]
and hence it is clear that $d_i=1$. Here $R_\phi$ denotes the plane rotation by angle $\phi$.
Comparing (\ref{restric}) and (\ref{restric1}) we obtain
\[f^*(\sigma_1) = g^* h^*(\sigma_1) = g^*(\sum_{i=1}^n (-1)^{i+1}s_i ) = s,\]
where we have used the assumption that $n$ is odd.
 (Note that for $n$ even the above arguments give $f^\ast(\sigma_1)=0$.)

Let us examine the Serre spectral sequence of fibration (\ref{fibration1}).
It has two rows and may have one nontrivial differential. Since we know that the fundamental class of the fiber
$s\in H^1(S^1;\kk)$ survives, i.e. application of the differential to it gives zero, it follows
that all the differentials in the Serre
spectral sequence vanish. We conclude that the cohomology algebra of the base
$H^\ast(G(S^2,n);\kk)$ is the factor of 
 $H^\ast(SO(3)\times G(S^2;A,A,n-1);\kk)$ with respect to the ideal generated by class $\sigma_1$.
Since $\sigma_{2i+1} = \sigma_1\sigma_{2i}$, we obtain that 
$H^\ast(G(S^2,n);\kk)$ has generators $w$ with $\deg w =3$ and $\sigma_{2i}$, with $\deg \sigma_{2i} =2i$,
where $i=0, 1,\dots, (n-3)/2$, which satisfy
relations (\ref{relnew}).

\end{proof}

\section{\bf Proof of Theorem \ref{thm5}}

Let $T\subset \R^{m+1}$ be a compact strictly convex domain with smooth boundary $X=\partial T$.
Consider the smooth function
\begin{eqnarray*}
L_X : G(X, n)\to \R, \quad L_X(x_1, \dots, x_n) = -\sum\limits_{i=1}^{i=n}|x_i -x_{i+1}|
\end{eqnarray*}
(the negative total length), where we understand the indices cyclically modulo $n$, i.e. $x_{n+1}=x_{1}$. 
The critical points of $L_X$ are in 1-1 correspondence with $n$-periodic billiard trajectories in $X$.

Fix  $\epsilon>0$
and consider 
$$G_\epsilon\subset G(X, n), \quad G_\epsilon =\{(x_1, \dots, x_n)\in X^{\times n}: \prod\limits_{i=1}^{n} |x_i-x_{i+1}| \ge \epsilon\}.
$$
According to Proposition 4.1 from \cite{FT}, if $\epsilon>0$ is small enough then: 

\begin{itemize}
\item[(a)] $G_\epsilon$ is a compact manifold with boundary;

\item[(b)]  the inclusion $G_\epsilon\subset G(X,n)$ is a $D_n$-equivariant homotopy equivalence;

\item[(c)]  all critical points of $L_X$ are contained in $G_\epsilon$;

\item[(d)]  at every point of $\partial G_\epsilon$ the gradient of $L_X$ has the outward direction.
\end{itemize}

Now we apply Proposition \ref{prop16} with $M=G_\epsilon$, $f=L_X$, and $G=D_n$.
We conclude that the number of $D_n$-orbits of $n$-periodic billiard trajectories in $X$ is at least 
\[\cat(G_\epsilon/D_n) = \cat(G(X,n)/D_n).\]
Since we assume that $n$ is a odd prime, the action of $D_n$ on $G(X,n)$ is free, and we may use inequality 
(\ref{cover}) which gives 
\[\cat(G(X,n)/D_n) \ge \cat(G(X,n))\ge \cl(G(S^m,n)) +1.\]

Theorems \ref{openstring6} and \ref{openstring8} allow to estimate $\cat(G(S^m,n)$. 
Assume first that $m$ is odd, $m>1$. Then (according to Theorem \ref{openstring6})
we have a nonzero cup-product 
$$\sigma_1^{n-2}u,$$ 
which shows that the cup-length
of $G(S^m,n)$ for odd $m>1$ is at least $n-1$. This gives lower bound $n$ on the number of $D_n$-orbits
of $n$-periodic billiard trajectories in $X$ for $m$ odd.

If $m$ is even, then (by Theorem \ref{openstring8}) we have a nonzero cup-product 
$$\sigma_2^{\frac{n-3}{2}}w\, 
\in H^\ast(G(S^m,n);\kk),$$ 
where $\kk$ is a field of characteristic $\ne 2$.
This shows that for even $m$ the cup-length
of $G(S^m,n)$ for even $m$ is at least $(n-1)/2$. This gives lower bound $(n+1)/2$ on the number of $D_n$-orbits
of $n$-periodic billiard trajectories in $X$ for $m$ even.

\qed

\bibliographystyle{amsalpha}

\end{document}